\documentclass{siamonline1116}

\usepackage{lipsum}
\usepackage{amsfonts}
\usepackage[mathscr]{euscript}
\usepackage{graphicx}
\usepackage{epstopdf}
\usepackage{tikz}
\usepackage{onimage}
\usepackage{algorithm}
\usepackage{algorithmic}
\usepackage{amsmath}
\usepackage{amssymb}
\usepackage{latexsym}
\usepackage{bbm}
\usepackage{relsize}
\usepackage[caption=false]{subfig}
\usepackage{multirow}
\usepackage{soul}
\ifpdf
  \DeclareGraphicsExtensions{.eps,.pdf,.png,.jpg}
\else
  \DeclareGraphicsExtensions{.eps}
\fi

\numberwithin{theorem}{section}

\newcommand{\TheTitle}{Learning Bilinear Models of Actuated Koopman Generators from Partially-Observed Trajectories}
\newcommand{\ShortTitle}{Learning Koopman Generators using the EM Algorithm}
\newcommand{\TheAuthors}{S. E. Otto, S. Peitz, and C. W. Rowley}

\headers{\ShortTitle}{\TheAuthors}

\title{{\TheTitle}\thanks{Last updated \today.}}

\author{
  Samuel E. Otto\thanks{Mechanical and Aerospace Engineering, Princeton University, Princeton, NJ
    (\email{sotto@princeton.com}).} 
    \and
    Sebastian Peitz\thanks{Department of Computer Science, Paderborn University, Germany
	(\email{sebastian.peitz@upb.de}).}
    \and
    Clarence W. Rowley \thanks{Mechanical and Aerospace Engineering, Princeton University, Princeton, NJ
    (\email{cwrowley@princeton.edu}).}
}

\usepackage{amsopn}

\DeclareMathOperator*{\minimize}{\min\!imize\enskip}

\DeclareMathOperator{\Tr}{Tr}

\DeclareMathOperator{\avg}{avg}
\DeclareMathOperator{\grad}{\nabla}
\DeclareMathOperator{\res}{res}

\newcommand{\mat}[1]{{\mathbf{#1}}}
\newcommand{\vect}[1]{{\mathbf{#1}}}
\newcommand{\vmu}{\boldsymbol{\mu}}
\newcommand{\vmuhat}{\boldsymbol{\hat{\mu}}}
\newcommand{\mSig}{\boldsymbol{\Sigma}}
\newcommand{\mSighat}{\boldsymbol{\hat{\Sigma}}}
\newcommand{\vpsi}{\boldsymbol{\psi}}

\newcommand{\vf}{\vect{f}}
\newcommand{\vg}{\vect{g}}

\newcommand{\vu}{\vect{u}}
\newcommand{\vv}{\vect{v}}
\newcommand{\vw}{\vect{w}}
\newcommand{\vx}{\vect{x}}
\newcommand{\vy}{\vect{y}}
\newcommand{\vz}{\vect{z}}
\newcommand{\calP}{\mathcal{P}}

\newcommand{\mQ}{\boldsymbol{Q}}
\newcommand{\mR}{\boldsymbol{R}}

\newsiamremark{example}{Example}
\newsiamremark{remark}{Remark}

\makeatletter
\renewcommand*\env@matrix[1][*\c@MaxMatrixCols c]{%
  \hskip -\arraycolsep
  \let\@ifnextchar\new@ifnextchar
  \array{#1}}
\makeatother

\makeatletter
\newcounter{phase}[algorithm]

\makeatother

\ifpdf
\hypersetup{
  pdftitle={\TheTitle},
  pdfauthor={\TheAuthors}
}
\fi

\begin{document}

\maketitle

\begin{abstract}
Data-driven models for nonlinear dynamical systems based on approximating the underlying Koopman operator or generator have proven to be successful tools for forecasting, feature learning, state estimation, and control.
It has become well known that the Koopman generators for control-affine systems also have affine dependence on the input, leading to convenient finite-dimensional bilinear approximations of the dynamics.
Yet there are still two main obstacles that limit the scope of current approaches for approximating the Koopman generators of systems with actuation.
First, the performance of existing methods depends heavily on the choice of basis functions over which the Koopman generator is to be approximated; and there is currently no universal way to choose them for systems that are not measure preserving.
Secondly, if we do not observe the full state, then it becomes necessary to account for the dependence of the output time series on the sequence of supplied inputs when constructing observables to approximate Koopman operators.
To address these issues, we write the dynamics of observables governed by the Koopman generator
as a bilinear hidden Markov model, and
determine the model parameters using the expectation-maximization (EM) algorithm.
The E-step involves a standard Kalman filter and smoother, while the M-step resembles control-affine dynamic mode decomposition for the generator.
We demonstrate the performance of this method on three examples, including recovery of a finite-dimensional Koopman-invariant subspace for an actuated system with a slow manifold; estimation of Koopman eigenfunctions for the unforced Duffing equation; and model-predictive control of a fluidic pinball system based only on noisy observations of lift and drag.
\end{abstract}



\section{Introduction}

Data-driven modeling techniques based on approximating the underlying Koopman operators of dynamical systems are promising methods for forecasting, estimation, and control.
Rather than directly predicting the nonlinear trajectories of states, the Koopman operator describes the evolution of functions on the state space; this evolution operator is always linear, but possibly infinite-dimensional \cite{Rowley:2009}.
In practice, one often makes a truncated approximation of the Koopman operator on a finite-dimensional space of functions, leading to a linear model for the dynamics of a lifted state composed of basis function values \cite{Schmid:2010,Williams2016extending}.
When accurate low-dimensional models of this type can be found, they enable efficient forecasting and state estimation as well as synthesis of control strategies that can be implemented in real time \cite{Proctor2016dynamic,Korda2018linear,Peitz2019koopman,Peitz2020data}.
This type of approach is especially attractive for high-dimensional systems like fluid flows and other spatiotemporal equations, where the models reflect the dynamics of a smaller collection of salient features that remain coherent in space and time.
For recent reviews of these topics, see \cite{Mauroy2020koopman, Otto2021koopman, Brunton2021modern}.

The most significant obstacle to producing an accurate low-dimensional model based on an approximation of the Koopman operator is the need for an appropriate finite-dimensional space of functions in which to make the approximation.
In particular, for the high-dimensional systems of interest one cannot rely on detailed analyses of the state space that may be performed in low dimensions to inform the choice of functions.
Instead, the prevailing approach is to rely on data collected from the system, and this is the approach we consider here as well.
Let us begin by mentioning the non-trivial cases when there are known to be good solutions to this problem.
When the underlying system is measure-preserving, the techniques introduced in \cite{Das2021reproducing, Govindarajan2019approximation, Colbrook2022mpEDMD} yield approximations of Koopman operators with guaranteed convergence properties for the spectrum and spectral projections.
These techniques enable short-horizon forecasts to be made even for chaotic dynamical systems.
However, the introduction of control or even the study of transients destroys the measure-preserving assumptions.
For autonomous, non-measure-preserving systems where the full state (or a diffeomorphic copy obtained via Takens embedding) can be observed, the ResDMD algorithm in \cite{Colbrook2021rigorous, Colbrook2023residual} provides similar convergence guarantees.
In practice, the choice of subspace (spanned by a dictionary of user-specified functions) in which to approximate Koopman operators, is still important.
Approaches based on deep learning \cite{Takeishi2017learning, Lusch2018deep, Wehmeyer2018time, Morton2018deep, Yeung2019learning, Otto2019linearly, Pan2020physics} have demonstrated the ability to identify subspaces of functions that are nearly invariant under the Koopman operator in the sense that
the the Koopman operator acting on these functions has a small component orthogonal to the subspace (in an $L^2$ sense).
These methods allow for accurate short-horizon predictions that are critical for engineering applications like estimation and control.

The situation when the full state is not available and we only have access to actuated trajectories of noisy partial observations of the state is significantly more difficult.
However, it is also the one most likely to be encountered in experimental applications where one forces the system with variety of input signals, but does not have access to enough sensors to measure the full state.
One approach could be to the use of delay coordinates, as in \cite{Kamb2020time, Das2019delay}, to obtain a diffeomorphic copy of the state via the Takens embedding theorem \cite{Takens1981detecting, Noakes1991takens}.
These ideas might be extended to actuated systems using non-autonomous delay embedding theorems \cite{Stark1999delay, Robinson2008topological}.
However, one must contend with the fact that the delayed observables now depend (in a complicated way) on the supplied input sequence in addition to the system's state.
When the equations of the original system are known, 
information about the state can be recovered from the delayed output sequence and the known input sequence via an appropriate state estimator \cite{Stengel1994optimal}.
When the governing equations are not available, our approach aims to construct a surrogate model of the system and to estimate the state variables of this model using data coming only from input-output sequences.

In this paper, we consider input-affine systems (defined in \cref{sec:bilinear}), for which it is known that the corresponding Koopman generators are also input-affine \cite{Surana2016koopman, Goswami2017global, Huang2020data, Peitz2020data,NPP+21}.
Approximating the input-affine Koopman generator in a finite-dimensional space of functions leads to a bilinear model of the dynamics.
A large class of systems can be approximated in this way, and there are even conditions under which an exact finite-dimensional bilinear model can be obtained \cite{Svoronos1980bilinear, Lo1975global}.
Even if the system has more complicated dependence on the input, the dynamics can be expanded as an affine combination of autonomous vector fields scaled by nonlinear functions of the input.
These nonlinear functions of the input can be treated as a new lifted input, for which the system is input-affine.
We note that it has become common to make Linear Time-Invariant (LTI) approximations of the lifted dynamics \cite{Proctor2016dynamic, Korda2018linear, Arbabi2018data, Narasingam2019koopman, Kaiser2020data} due to the abundance of available system identification algorithms like the Eigensystem Realization Algorithm (ERA) \cite{Juang1985eigensystem} and Observer/Kalman filter Identification (OKID) \cite{Juang1993identification} as well as powerful techniques for control synthesis \cite{Skogestad2007multivariable}.
However, these ``lifted LTI'' systems must obey the rather restrictive principle of linear superposition, regardless of their dimension.
Given this fact, we provide an example (\cref{ex:ltibad}) that illustrates why it is often necessary to consider the larger class of bilinear systems.
The advantages of bilinear models based on the Koopman generator for control applications are discussed in \cite{Bruder2021advantages}.

To remedy the problems associated with input dependence and partial observations in Koopman generator approximations, we treat the observed trajectories as if they come from a hidden Markov model (HMM) whose hidden dynamics are governed by an input-affine combination of finite-dimensional Koopman generators acting on the same space.
The matrix approximations of these Koopman generators, along with the error covariances of the model, are treated as parameters which are learned from the provided trajectory data using an Expectation-Maximization (EM) algorithm \cite{Baum1970maximization, Shumway1982approach, Ghahramani1996parameter}.
By taking this approach, we avoid specifying the shared space of functions ahead of time.
Rather, these functions are implicitly defined by the Koopman generator approximations, and their values are estimated during the expectation (E) step of the EM algorithm using an optimal state estimator.
The estimator can be viewed as untangling the complicated dependence of the observed output sequences on the initial states and the time-varying input.
The maximization (M) step of the EM algorithm uses these untangled state estimates to construct updated matrix approximations of the Koopman generators using a procedure that very closely resembles the Extended Dynamic Mode Decomposition (EDMD)-based technique for input-affine Koopman generators described in Peitz et al. \cite{Peitz2020data}.
While this procedure appears to be intuitive, the fact that it can be derived from the EM algorithm guarantees that the likelihood function of the model increases at every iteration.

Related work includes bilinear system identification algorithms including an EM algorithm \cite{Gibson2005maximum} and a variant of Observer/Kalman filter Identification (OKID) \cite{Vicario2014okid}.
While our approach is similar to the one in \cite{Gibson2005maximum}, our model inherits special structure from the Koopman generator framework.
In particular, we include a state that remains constant for all time in addition to simplifying how the outputs are reconstructed from the hidden state in our model.
Moreover, we provide new connections between the EM algorithm, EDMD, and approximation techniques for the Koopman generator.

The plan for our paper is as follows:
In \cref{sec:bilinear} we discuss the construction and advantages of bilinear surrogate models based on approximating the Koopman generators of input-affine dynamical systems.
\Cref{sec:formulating_the_HMM} presents the formulation of the Hidden Markov Model we use to approximate Koopman generators.
In \cref{sec:EM_algorithm} we describe the Expectation-Maximization (EM) algorithm that we use to learn the parameters of the model from input-output data along trajectories.
We also discuss how the Expectation (E) and Maximization (M) steps are related to Extended Dynamic Model Decomposition and Kalman filtering.
In \cref{sec:approximating_eigenfunction} we discuss how the filtering and smoothing procedure in the E step can be used to approximate Koopman eigenfunctions.
Finally, in \cref{sec:numerical_applications}, we present several numerical applications, including approximation of Koopman spectra and model-predictive control.

\section{Bilinear models based on approximating the Koopman generator}
\label{sec:bilinear}

Consider a dynamical system on a state space $\mathcal{X}$ whose time-$t$ flow map under actuation~$u$ is given by $\Phi_u^t:\mathcal{X} \to \mathcal{X}$.
For the time being, we take the actuation~$u$ to be constant; we will allow time-varying inputs below, once we introduce the Koopman generator.
If $f: \mathcal{X} \to \mathbb{C}$ is a complex-valued function on the state space, then we can define another function on the state space by the composition $f \circ \Phi_u^t$.
Over any function space $\mathcal{F}$ that is closed under this composition operation, we can define the Koopman operator according to
\begin{equation}
    U^t_u f := f \circ \Phi_u^t.
\end{equation}
While the original state space $\mathcal{X}$ need not have any structure and the dynamics $\Phi_u^t$ may be very complicated, the Koopman operator is always a linear operator, albeit on a possibly infinite-dimensional space $\mathcal{F}$.
By taking the perspective of functions on the state space, sometimes called ``observables,'' the Koopman operator allows us to study the aggregate, global behavior of the system, rather than individual trajectories.
Of particular interest is the spectrum of the Koopman operator, which can be used as an indicator of periodic or mixing dynamics, in addition to any Koopman eigenfunctions which tell us a great deal about the state space geometry, including revealing attractors \cite{Mezic:2019, Mauroy2016global, Mauroy2013isostables}.
Moreover, if the elements of $\vpsi = (\psi_1, \ldots, \psi_N)$, $\psi_j\in\mathcal{F}$ span a subspace that is invariant under the Koopman operator, then the evolution of these quantities is governed by a linear system; i.e., there is a matrix $\mat{U}^t_u \in \mathbb{C}^{N\times N}$ such that
\begin{equation}
    U^t_u\vpsi := \vpsi \circ \Phi_u^t = \left(\mat{U}^t_u\right)^T \vpsi.
\end{equation}
Even when there are no non-trivial finite-dimensional invariant subspaces (as in the case of mixing dynamics),
finite-dimensional approximations of $U^t_u$ can be used to rigorously approximate the spectrum, spectral measures, and functional calculus of Koopman operators \cite{Das2021reproducing, Govindarajan2019approximation, Colbrook2021rigorous, Colbrook2022mpEDMD}.
In particular, eigenfunctions of finite-dimensional approximations often provide useful observables for marking short-time forecasts of the system.

If the span of the elements of $\vpsi$ does not form an invariant subspace of the Koopman operator, there are a variety of techniques for approximating their evolution under the Koopman operator back in their span.
For instance, Extended Dynamic Mode Decomposition (EDMD) \cite{Williams2015data} is a widely used technique that can be viewed as a data-driven Petrov-Galerkin projection of the Koopman operator onto the span of a collection of user-defined functions.
However, in the actuated setting, EDMD would require computing a new matrix $\mat{U}_u^t$ for every input level $u$ under consideration.
A related data-driven technique developed by \cite{Williams2016extending} overcomes this problem by postulating a form for the dependence of $\mat{U}_u^t$ on $u$ and solving for the unknown coefficients to best fit the data in a least-squares sense.

Disadvantages of working with a family of Koopman operators $U_u^t$ parameterized by the input $u$ include the complicated dependence of $U_u^t$ on the input and the fact that the input must be held constant over the time interval.
Instead, we may consider the Koopman generator $V_u : D(V_u) \subset \mathcal{F} \to \mathcal{F}$, which is a linear operator yielding the time derivative of functions $f$ in its domain according to
\begin{equation}
    V_u f := \left. \frac{\mathrm{d}}{\mathrm{d}t} U_u^t f \right\vert_{t=0}.
\end{equation}
When the underlying vector field that generates the flow $\Phi_u^t$ is input-affine, that is,
\begin{equation}\label{eq:bilinear_system}
    \dot{x} = F_0(x) + \sum_{k=1}^{\dim u} u_k(t) F_i(x),
\end{equation}
then the Koopman generator applied to a differentiable function $f$ is given by the Lie derivative
\begin{equation}
    V_u f = \underbrace{F_0 \cdot \grad}_{V_0} f + \sum_{k=1}^{\dim u} u_k(t) \underbrace{F_k \cdot \grad}_{V_k} f.
\end{equation}
In such a situation, the dynamics of differentiable functions $f$ are bilinear,
where each operator $V_k = F_k \cdot \grad$ is the Koopman generator of its respective component vector field $F_k$.

There are a variety of techniques for approximating the dynamics of a collection of observables $\vpsi = (\psi_1, \ldots, \psi_N)$ under the action of the Koopman generator back in the span of these observables.
These techniques, which include Carleman linearization \cite{Carleman1932application, Kowalski1991nonlinear, Brunton2016koopman} and EDMD techniques for the generator \cite{Klus2020data, Peitz2020data}, seek to identify the $N\times N$ matrices $\mat{V}_0, \mat{V}_1, \ldots, \mat{V}_{\dim u}$ in a bilinear approximation
\begin{equation}
    V_u \vpsi \approx \mat{V}_0^T \vpsi + \sum_{k=1}^{\dim u} u_k \mat{V}_k^T \vpsi.
\end{equation}

In many cases, we are only interested in the dynamics of some relevant observed quantities
\begin{equation}
    \vect{y} = \vect{g}(x),
\end{equation}
where $\vect{g}$ might not be injective, thus revealing incomplete information about the underlying state $x$.
It is natural to approximate the observed quantities in the span of our dictionary $\vpsi$ according to $\vect{g} = \mat{C} \vpsi$ for some matrix $\mat{C}$, which we may also seek to identify, e.g, by using Petrov-Galerkin projection.
The situation when $\vect{g}$ is not injective and our data consists only of these observed quantities $\{ \vect{y}_1, \ldots \vect{y}_L\}$ along actuated trajectories poses serious problems for existing identification techniques based on approximating the underlying Koopman operator or generator as we shall see below in \cref{subsec:challenge_with_partial_observations}.
In particular, we seek to identify matrices $\{\mat{V}_k\}$ and $\mat{C}$ along with the functions $\vpsi$ so that the model
\begin{equation}
\boxed{
\begin{aligned}
    \vect{\dot{z}} &= \mat{V}_0^T\vect{z} + \sum_{k=1}^{\dim u} u_k \mat{V}_k^T\vect{z}, \qquad \vect{z} = \vpsi(x) \\
    \vect{y} &= \mat{C} \vect{z}
\end{aligned}
}
\label{eqn:finite_dimensional_generator_IO_model}
\end{equation}
is in close agreement with our observed data $\{ \vect{y}_1, \ldots \vect{y}_L\}$ and most importantly is capable of predicting the dynamics of $\vect{y}(t)$ along new trajectories given its history.

\subsection{Limitations of lifted linear time-invariant models}
A special case of the lifted bilinear model \cref{eqn:finite_dimensional_generator_IO_model} that has received a lot of attention in the literature \cite{Proctor2016dynamic, Korda2018linear, Arbabi2018data, Narasingam2019koopman, Kaiser2020data} is the lifted linear-time-invariant (LTI) model
\begin{equation}
\begin{aligned}
    \vect{\dot{z}} &= \mat{A} \vect{z} + \mat{B}\vect{u}, \qquad \vect{z} = \vpsi(x) \\
    \vect{y} &= \mat{C}\vect{z}.
\end{aligned}
\label{eqn:lifted_LTI_model}
\end{equation}
To see that this is a special case of the bilinear model \cref{eqn:finite_dimensional_generator_IO_model}, we re-write \cref{eqn:lifted_LTI_model} using observables $(1, \psi_1, \ldots, \psi_N)$ that evolve according to
\begin{equation}
    \begin{bmatrix}
    0 \\
    \vect{\dot{z}}
    \end{bmatrix} = 
    \begin{bmatrix}
    0 & \vect{0} \\
    \vect{0} & \mat{A}
    \end{bmatrix}
    \begin{bmatrix}
    1 \\
    \vect{z}
    \end{bmatrix}
    + \sum_{k=1}^{\dim u} u_k
    \begin{bmatrix}
    0 & \vect{0} \\
    \vect{b}_k & \mat{0}
    \end{bmatrix}
    \begin{bmatrix}
    1 \\
    \vect{z}
    \end{bmatrix}, \qquad
    \mat{B} = \begin{bmatrix}
    \vect{b}_1 & \cdots & \vect{b}_{\dim u}
    \end{bmatrix}.
    \label{eqn:lifted_LTI_model_as_bilinear}
\end{equation}
Lifted LTI models are appealing because standard machinery from linear systems and control theory can be applied to them directly.
This includes system identification techniques like the Eigensystem Realization Algorithm (ERA) \cite{Juang1985eigensystem} and Observer/Kalman Filter Identification (OKID) \cite{Juang1993identification} that allow $\mat{A}$ and $\mat{B}$ to be determined from trajectories of the observed output.

While lifted LTI models are a strictly broader class than LTI systems on the original state space, the dynamics of observed quantities $\vect{y}$ must still obey the principle of linear superposition.
This is a significant limitation as the following example demonstrates.
\begin{example}[A system that doesn't admit an accurate lifted LTI model]
\label{ex:ltibad}
There are some very simple systems, for instance
\begin{equation*}
\begin{aligned}
    \dot{x} &= u x, \qquad x(0) = 1 \\
    y &= x ,
\end{aligned}
\end{equation*}
that violate the superposition principle so badly that there cannot exist an accurate lifted LTI model in the form of \cref{eqn:lifted_LTI_model}, regardless of how large the dimension, $\dim \vect{z}$, of such an approximation is taken to be.
To see why this system does not admit a lifted LTI approximation, consider the three trajectories
\begin{align*}
    u_a \equiv -1 \quad &\Rightarrow \quad x_a(t) = e^{-t} \\
    u_b \equiv -3 \quad &\Rightarrow \quad x_b(t) = e^{-3t} \\
    u_c \equiv 1 \quad &\Rightarrow \quad x_c(t) = e^{t}.
\end{align*}
and observe that for any dictionary of observables $\vpsi$, we have
\begin{align*}
    \vpsi(x_c(0)) &= 2 \vpsi(x_a(0)) - \vpsi(x_b(0)) = \vpsi(1) \\
    u_c &= 2 u_a - u_b.
\end{align*}
Therefore, if a lifted LTI model in the form of \cref{eqn:lifted_LTI_model} agrees with the first two trajectories $x_a(t)$ and $x_b(t)$, then it must predict the third trajectory to be
\begin{equation*}
    x_c^{(\text{LTI})}(t) = 2 x_a(t) - x_b(t) = 2 e^{-t} - e^{-3t}
\end{equation*}
by linear superposition.
This is a very poor prediction because it decays to zero exponentially, 
whereas, the real trajectory $x_c(t)$ blows up exponentially.
\end{example}

\subsection{Challenges with identifying partially-observed actuated systems}
\label{subsec:challenge_with_partial_observations}
One common feature of the existing techniques for approximating the Koopman operator and generator is their dependence on the choice of dictionary functions making up $\vpsi$.
This choice strongly affects the performance of these algorithms both in terms of predictive accuracy and approximating the spectrum of the Koopman operator.
If the elements of $\vpsi$ span a Koopman-invariant subspace, then any eigenvalue of the matrix $\mat{U}_u^t$ identified by EDMD is also an eigenvalue of the Koopman operator corresponding to an eigenfunction of $U_u^t$ in the span of $\psi_1,\ldots, \psi_N$.
Generically, however, the matrix $\mat{U}_u^t$ that we compute numerically always has $N$ distinct eigenvalues, but these need not correspond to eigenvalues of the Koopman operator if $\psi_1,\ldots, \psi_N$ do not span an invariant subspace.
In fact, when the underlying system exhibits a type of chaos knows as weak mixing \cite{Cornfeld:1982}, the Koopman operator has only a single eigenfunction, equal to a constant, and corresponding to eigenvalue $1$.

Recognizing the importance of the choice of the dictionary $\vpsi$, techniques have been developed to learn appropriate functions from data by parameterizing them using neural networks \cite{Takeishi2017learning, Lusch2018deep, Yeung2019learning, Otto2019linearly, Pan2020physics}.
Yet even these state-of-the-art methods have significant limitations: in particular, they require large amounts of training data consisting of rich features that are sufficient to determine the full state of the system.

Identifying the parameters needed to build a model in the form of \cref{eqn:finite_dimensional_generator_IO_model} is especially difficult because each observation $\vect{y}_l$ does not contain enough information to reconstruct the state $x$.
This constrains the dictionary elements we can construct to those that are expressible as $\psi_i(x) = \tilde{\psi}_i(\vect{g}(x))$ for some $\tilde{\psi}_i$.
When there is no actuation, the situation is easily solved by employing Takens' embedding theorem \cite{Takens1981detecting, Noakes1991takens} to construct functions of time-delayed observables
\begin{equation}
    \psi_i = \tilde{\psi}_i\left( \vect{g},\ \vect{g}\circ \Phi^{\Delta t},\ \ldots,\ \vect{g}\circ \Phi^{m \Delta t} \right),
    \label{eqn:delay_observables}
\end{equation}
as is done in \cite{Kamb2020time} and \cite{Das2019delay}.
With sufficiently many time-delays, $m$, and under certain genericity assumptions (see \cite{Takens1981detecting, Noakes1991takens}) the map 
$$x\mapsto \left( \vect{g}(x),\ \vect{g}(\Phi^{\Delta t}(x)),\ \ldots,\ \vect{g}(\Phi^{m \Delta t}(x)) \right)$$
becomes injective and functions on the state space can be written in the form of \cref{eqn:delay_observables}.
Most importantly, this allows one to search over functions on the state space implicitly by considering functions $\tilde{\psi}_i(\vect{y}(t_0),\ \vect{y}(t + \Delta t),\ \ldots,\ \vect{y}(t_0 + m \Delta t))$ formed from consecutive observations that one might have access to from an experiment.

However, for actuated systems the time-delayed output sequences now also depend on the input signal applied to the system.
For an actuated system with piecewise-constant inputs $u_1, \ldots, u_m$ along each time interval, the time-delayed observables take the form
\begin{equation}
    \psi_i = \tilde{\psi}_i\left( \vect{g},\ \vect{g} \circ \Phi_{u_1}^{\Delta t},\ \vect{g} \circ \Phi_{u_2}^{\Delta t} \circ \Phi_{u_1}^{\Delta t} ,\ \ldots,\ \vect{g} \circ \Phi_{u_m}^{\Delta t} \circ \cdots \circ \Phi_{u_2}^{\Delta t} \circ \Phi_{u_1}^{\Delta t} \right).
    \label{eqn:delay_observables_actuation}
\end{equation}
Delay embedding theorems for non-autonomous systems ensure that under certain genericity assumptions, a sufficiently long delayed output sequence provides an input-dependent embedding of the system's state \cite{Robinson2008topological, Stark1999delay}.
However, we must contend with the fact that the observable defined by \cref{eqn:delay_observables_actuation} depends on both the state and the supplied input sequence.
On the one hand, we may not know how well these particular observables work in a specific situation. 
On the other hand, there may be good reasons for a specific choice of other observables; motivated, e.g., by expert knowledge.

An approach described in \cite{Korda2018linear} seeks to approximate the Koopman operator of the autonomous system describing the evolution of an augmented state variable that combines the system's state with the infinite sequence of inputs applied over time.
Certain assumptions, e.g., concerning the form of the observables, are needed to make the approximation problem for the augmented Koopman operator tractable and useful for control applications in which the entire input sequence is not known ahead of time.

We take a different approach that seeks to untangle the effects of actuation and the initial condition on observed time histories.
As we will see in the sequel, a standard Kalman filter and smoother is the appropriate tool for accomplishing this task when the model is bilinear as in \cref{eqn:finite_dimensional_generator_IO_model}.
This state estimator will be used as the first step in an iterative two-step algorithm for learning the parameters of the model \cref{eqn:finite_dimensional_generator_IO_model}.
The second step resembles the actuated EDMD method for the generator presented in \cite{Peitz2020data}.

\section{Koopman generator approximations as parameters in a hidden Markov model}
\label{sec:formulating_the_HMM}
Let us suppose that we are given a collection of finely sampled time histories of inputs $\{ \vect{u}_l^{(m)} \}_{l=0}^{L-1}$ and observations $\{ \vect{y}_l^{(m)} \}_{l=0}^L$ along trajectories $m=1,\ldots,M$ from a control-affine system.
We seek to identify matrix approximations of the Koopman generators $\{ \mat{V}_k\}$ together with an appropriate collection of functions $\vpsi$ so that a bilinear model in the form of \cref{eqn:finite_dimensional_generator_IO_model} fits the data as closely as possible (in the sense of maximum likelihood described later in \cref{sec:EM_algorithm}).
While we do not know $\vpsi$ ahead of time, we will always force the first function $\psi_1$ to be constant since the constant function is always an eigenfunction of the Koopman operator.
In the setting of a bilinear model \cref{eqn:finite_dimensional_generator_IO_model}, forcing one of the states to be constant allows for linear input terms to be modeled as in \cref{eqn:lifted_LTI_model_as_bilinear}, as well as allowing for a constant shift to be applied to the output.
If the observations $\vect{g}$ are linearly independent and expressible in the
span of the yet unknown $\vpsi = (1, \boldsymbol{\tilde{\psi}})$, we can always
find a linear transformation of $\boldsymbol{\tilde{\psi}}$ such that
\begin{equation}
    \mat{C} = \begin{bmatrix}
    \vect{c}_0 & \mat{I}_{\dim\vect{y}} & \mat{0}
    \end{bmatrix},
\end{equation}
where only the column $\vect{c}_0$ is unknown.
Therefore, we will constrain the matrix $\mat{C}$ to always have the above form.
For ease of notation, we will also assume that the first component of the input satisfies $u_0 \equiv 1$ so that we do not have to treat the drift term separately.

Since the time sampling of our observations is assumed to be fine, we propose to work with a hidden Markov model obtained by explicit Euler discretization of \cref{eqn:finite_dimensional_generator_IO_model}.
For simplicity, we assume that the sampling time interval $\Delta t$ is fixed.
Furthermore, we assume that the combined process noise and modeling error can be modeled using independent Gaussian random vectors $\vect{w}_l^{(m)} \sim \mathcal{N}(\vect{0}, \mSig_{\vect{w}})$ and the combined measurement noise and modeling error can be modeled using independent Gaussian random vectors $\vect{v}_l^{(m)} \sim \mathcal{N}(\vect{0}, \mSig_{\vect{v}})$.
These assumptions yield the bilinear HMM
\begin{equation}
\boxed{
\begin{aligned}
    \begin{bmatrix}
    1 \\
    \vect{z}_{l+1}
    \end{bmatrix}
    &= \left( \mat{I} + \Delta t \sum_{k=0}^{\dim \vect{u}} [\vect{u}_l]_k \mat{V}_k \right)^T
    \begin{bmatrix}
    1 \\
    \vect{z}_{l}
    \end{bmatrix}
    + \begin{bmatrix}
    0 \\
    \vect{w}_{l}
    \end{bmatrix} \\
    \vect{y}_{l} &= 
    \vect{c}_0 + \begin{bmatrix} \mat{I} & \mat{0} \end{bmatrix} \vect{z}_l + \vect{v}_l
\end{aligned}
}
\label{eqn:main_HMM}
\end{equation}
where the matrix approximations of the Koopman generators are constrained to have their first column equal to zero, that is, $\mat{V}_k = \begin{bmatrix} \vect{0} & \mat{\tilde{V}}_k \end{bmatrix}$.
Another convenient way to express the latent state dynamics of our HMM \cref{eqn:main_HMM} is
\begin{equation}
    \vect{z}_{l+1} = \mat{A}_l \vect{z}_l + \vect{b}_l + \vect{w}_l,
\end{equation}
where $\mat{A}_l$ and $\vect{b}_l$ are determined at each time step by the corresponding entries of the first-order matrix approximation of the Koopman operator over the $l$th time interval
\begin{equation}
    \begin{bmatrix}
    1 & \vect{b}_l^T \\
    \vect{0} & \mat{A}_l^T
    \end{bmatrix} 
    = \mat{I} + \Delta t \sum_{k=0}^{\dim \vect{u}} [\vect{u}_l]_k \mat{V}_k
    =: \mat{U}_l.
    \label{eqn:finite_time_Koopman_operator_approx}
\end{equation}

\begin{remark}[Formulating the HMM using linear multi-step methods]
    \label{rem:linear_mulistep_methods}
    Our choice to use explicit Euler time discretization when formulating the HMM \cref{eqn:main_HMM} is partially motivated by the resulting simplicity of the EM algorithm and its connections to other common methods.
    However, it can introduce significant error for larger values of $\Delta t$.
    Using the Euler discretization, we will see that the maximization (M) step has an explicit solution resembling Extended Dynamic Mode Decomposition (EDMD) for the Koopman generator (see \cref{subsec:Mstep}) and the expectation (E) step is solved by an out-of-the-box Kalman filter and smoother (see \cref{subsec:Estep}).
    An explicit solution during the M-step, albeit a more complicated one, can still be obtained when \cref{eqn:main_HMM} is formulated using a linear multistep method of the form
    \begin{equation}
        \sum_{s=0}^{n_s} a_s \begin{bmatrix} 1 \\ \vect{z}_{l+s} \end{bmatrix}
        = \sum_{s=0}^{n_s} b_s \left( \sum_{k=0}^{\dim \vect{u}} [\vect{u}_{l+s}]_k \mat{V}_k \right)^T \begin{bmatrix} 1 \\ \vect{z}_{l+s} \end{bmatrix}
        + \begin{bmatrix} 0 \\ \vect{w}_{l} \end{bmatrix},
    \end{equation}
    determined by coefficients $\{a_s\}_{s=0}^{n_s}$ and $\{b_s\}_{s=0}^{n_s}$.
    While one can derive appropriate filtering equations for the E-step, the standard Kalman filter and smoother can no longer be applied directly.
    Using linear multistep methods may allow for improved approximations with larger time steps $\Delta t$.
    A systematic investigation is a subject for future work.
\end{remark}

The parameters of the model \cref{eqn:main_HMM} that we seek to determine are the generator matrices $\{ \mat{\tilde{V}}_k \}_{k=0}^{\dim\vect{u}}$, the observation offset vector $\vect{c}_0$, the noise covariance matrices $\mSig_{\vect{w}}$ and $\mSig_{\vect{v}}$, and a prior Gaussian distribution for the initial condition $\vect{z}_0 \sim \mathcal{N}(\vmu_0, \mSig_0)$; that is,
\begin{equation}
    \calP = \left( \mat{\tilde{V}}_0, \ldots, \mat{\tilde{V}}_{\dim\vect{u}},\ \vect{c}_0,\ \mSig_{\vect{w}},\ \mSig_{\vect{v}},\ \vmu_0, \mSig_0 \right).
\end{equation}
In order to identify these parameters, we will maximize the likelihood of generating our data $\{ \vect{y}_l^{(m)} \}$ under the assumed model \cref{eqn:main_HMM} using the expectation-maximization (EM) algorithm whose details are given in Section~\ref{sec:EM_algorithm}.

\section{EM algorithm for learning Koopman generator approximations}
\label{sec:EM_algorithm}

Maximum likelihood estimation entails maximizing the probability of the observed data over the model parameters.
In particular, let us combine all of our observations $\lbrace \vect{y}_l^{(m)} \rbrace_{0\leq l\leq L}$ along the $m$th independent trajectory into a matrix $\mat{Y}^{(m)}$ and denote the joint probability density of observing outputs $\mat{\tilde Y}$ under the model parameters $\calP$ by $P_{\mat{Y}}(\mat{\tilde{Y}} ;\ \calP)$.
We aim to maximize the log ``likelihood'' of the observations we collect, given by
\begin{equation}
    L(\calP) = \sum_{m=1}^M\log P_{\mat{Y}}(\mat{Y}^{(m)};\ \calP).
    \label{eqn:log_lik}
\end{equation}
The log likelihood is used instead of the raw likelihood because the probability $P_{\mat{Y}}(\mat{Y}^{(m)};\ \calP)$ factors into a product of many terms due to the Markov property of the model \cref{eqn:main_HMM}.

The required density $P_{\mat{Y}}(\mat{Y};\ \calP)$ can be expressed using the model \cref{eqn:main_HMM} by recognizing that it is a marginal distribution
\begin{equation}
    P_{\mat{Y}}(\mat{Y}^{(m)};\ \calP) = \int P_{\mat{Z}, \mat{Y}}(\mat{Z}, \mat{Y}^{(m)};\ \calP) d\mat{Z},
\end{equation}
where $P_{\mat{Z}, \mat{Y}}$ is the joint density of the observed trajectory $\mat{Y}$ and the hidden variables $\lbrace \vect{z}_l \rbrace_{l=0}^L$ stacked into a matrix $\mat{Z}$.
Yet this high-dimensional integral is rather difficult to evaluate and makes direct optimization of \cref{eqn:log_lik} futile.
By introducing a new probability density $\mat{Z} \mapsto Q^{(m)}(\mat{Z})$ that is to be determined and a random variable $\mat{\hat{Z}}^{(m)}$ with density $Q^{(m)}$, the above integral is converted into an expectation
\begin{equation}
    P_{\mat{Y}}(\mat{Y}^{(m)};\ \calP) = \mathbb{E}_{\mat{\hat{Z}}^{(m)}}\left[ \frac{P_{\mat{Z}, \mat{Y}}(\mat{\hat{Z}}^{(m)}, \mat{Y}^{(m)};\ \calP)}{Q(\mat{\hat{Z}}^{(m)})} \right].
\end{equation}
Recognizing that $\log$ is a concave function, we apply Jensen's inequality to obtain a lower bound on the log likelihood,
\begin{equation}
\boxed{
\begin{aligned}
    L(\calP) \geq \hat{L}_{Q}(\calP) &= \sum_{m=1}^M \left\{ \mathbb{E}_{\mat{\hat{Z}}^{(m)}}\left[ \log{P_{\mat{Z}, \mat{Y}}(\mat{\hat{Z}}^{(m)}, \mat{Y}^{(m)};\ \calP)} \right] - \mathbb{E}_{\mat{\hat{Z}}^{(m)}}\left[ \log{Q^{(m)}(\mat{\hat{Z}}^{(m)})} \right] \right\} \\
    &= L(\calP) - \sum_{m=1}^M D_{KL}\big(Q^{(m)} \Vert P_{\mat{Z}\vert \mat{Y}=\mat{Y}^{(m)}}\big),
\end{aligned}
}
\label{eqn:ELBO}
\end{equation}
which is commonly referred to as the variational lower bound \cite{Bishop2006pattern} or Evidence Lower Bound (ELBO) \cite{Blei2017variational}.
The quantity 
\begin{equation}
    D_{KL}\big(Q^{(m)} \Vert P_{\mat{Z}\vert \mat{Y}=\mat{Y}^{(m)}}\big) = \mathbb{E}_{\mat{\hat{Z}}^{(m)}}\left[ \log{\left(\frac{Q^{(m)}(\mat{\hat{Z}}^{(m)})}{P_{\mat{Z}\vert \mat{Y}=\mat{Y}^{(m)}}(\mat{\hat{Z}}^{(m)};\ \calP)}\right)} \right]
    \label{eqn:KL_div}
\end{equation}
is called the Kullback-Leibler (KL)-divergence or the entropy of $P_{\mat{Z}\vert \mat{Y}=\mat{Y}^{(m)}}$ relative to $Q^{(m)}$.

A key observation is that the inequality in \cref{eqn:ELBO} becomes equality when the probability density $Q^{(m)}$ is chosen to be the conditional density of $\mat{Z}$ given $\mat{Y}^{(m)}$; that is,
\begin{equation}
    Q^{(m)} = P_{\mat{Z}\vert \mat{Y}=\mat{Y}^{(m)}},\ m=1,\ldots, M \quad \Rightarrow \quad
    L(\calP) = \hat{L}_{Q}(\calP).
    \label{eqn:exact_inference}
\end{equation}
This is an immediate consequence of the definition of the KL-divergence in \cref{eqn:KL_div}.
For this reason, $Q$ is often referred to as the ``inference'' distribution since its optimal form allows one to infer the values of the latent variables from the observations.

The form of the ELBO \cref{eqn:ELBO} and the property \cref{eqn:exact_inference} suggest a kind of coordinate ascent optimization procedure for maximizing the likelihood that iteratively updates the inference distributions $Q^{(m)}$ with the parameters $\calP$ fixed and then updates the parameters $\calP$ with the inference distributions $Q^{(m)}$ fixed.
This is called the Expectation-Maximization (EM) algorithm \cite{Dempster1977maximum} because in the first step one computes the inference distribution at the current model parameters $\calP_0$ and uses it to assemble the ELBO $\calP \mapsto \hat{L}_Q(\calP)$ as an expectation \cref{eqn:ELBO} that is maximized in the second step.
This iterative method has the remarkable property that \emph{the log likelihood is guaranteed to increase at each step}: that is, if $\calP_1$ are the new model parameters found by the maximization step, then we have
\begin{equation}
    L(\calP_1) \geq \hat{L}_Q(\calP_1) \geq \hat{L}_Q(\calP_0) = L(\calP_0),
\end{equation}
where the first inequality is \cref{eqn:ELBO}, the second inequality is due to maximization, and the equality on the right is by definition of the expectation step and \cref{eqn:exact_inference}.
The convergence properties of the EM algorithm were studied by C. F. Jeff Wu in \cite{Wu1983convergence}, where it is shown that the limit points of the resulting sequence of iterates are stationary points of the likelihood function, assuming its superlevel set at the initial value is compact.
It is important to note that, because $\{ Q^{(m)} \}_{m=1}^M$ and $\calP$ are not jointly optimized, the EM algorithm only finds a local maximum of the likelihood when $\hat{L}_Q$ is a non-convex function of $\{ Q^{(m)} \}_{m=1}^M$ and $\calP$.
Unfortunately, this is a limitation in our case.
For more details about the EM algorithm, one can consult \cite{Bishop2006pattern} and~\cite{Blei2017variational}.

\subsection{Maximization step and connections with Dynamic Mode Decomposition}
\label{subsec:Mstep}
A key feature of the evidence lower bound given by \cref{eqn:ELBO} for our model \cref{eqn:main_HMM} is that maximization over the parameters $\mathcal{P}$ has an explicit solution.
Specifically, we maximize a regularized log likelihood function and evidence lower bound
\begin{equation}
    L(\calP) - R(\calP) \geq \hat{L}_Q(\calP) - R(\calP),
\end{equation}
where $L$ is given by \cref{eqn:log_lik}, $\hat{L}_Q$ is given by \cref{eqn:ELBO}, and the regularization function is defined by
\begin{equation}
    R(\calP) = 
     \frac{\gamma_{\mat{G}}}{2} \Delta t^2 \sum_{i=0}^{\dim\vu} \Tr\big( \mSig_{\vw}^{-1} \mat{\tilde{V}}_i^T \mat{\tilde{V}}_i \big)
    +  \frac{\gamma_{\vw}}{2} \Tr\big( \mSig_{\vw}^{-1} \big)
    +  \frac{\gamma_0}{2} \Tr\big( \mSig_{0}^{-1} \big)
    +  \frac{\gamma_{\vv}}{2} \Tr\big( \mSig_{\vv}^{-1} \big)
    \label{eqn:regularization_function}
\end{equation}
with nonnegative coefficients $\gamma_{\mat{G}}, \gamma_{\vw}, \gamma_0, \gamma_{\vv} \geq 0$.
We produce a maximizer of the regularized ELBO $\calP \mapsto \hat{L}_Q(\calP) - R(\calP)$ under a collection of positive-definiteness conditions that are guaranteed to hold when the coefficients in the regularization are taken to be positive.

Let the mean and joint covariance of the fixed inference distributions $Q^{(m)}$, $m=1,\ldots,M$ be denoted
\begin{equation}\label{eqn:muhatk_sighatk}
    \vmuhat_k^{(m)} = \mathbb{E}_{\mat{\hat{Z}}^{(m)}}\left[ \vect{\hat{z}}_k^{(m)} \right] 
    \qquad \text{and} \qquad
    \mSighat_{k,l}^{(m)} = \mathbb{E}_{\mat{\hat{Z}}^{(m)}}\left[ \left(\vect{\hat{z}}_k^{(m)} - \vmuhat_k^{(m)}\right) \left(\vect{\hat{z}}_l^{(m)} - \vmuhat_{l}^{(m)}\right)^T \right].
\end{equation}
With $\otimes$ denoting the Kronecker product, let
\begin{equation}
    \mat{G} = \sum_{m=1}^{M} \sum_{l=0}^{L-1} 
    \vu_l^{(m)} \otimes \big(\vu_l^{(m)}\big)^T \otimes
    \begin{bmatrix}
    1 & (\vmuhat_l^{(m)})^T \\
    \vmuhat_l^{(m)} & \mSighat_{l,l}^{(m)} + \vmuhat_l^{(m)}(\vmuhat_l^{(m)})^T 
    \end{bmatrix}
    + \gamma_{\mat{G}} \mat{I},
\end{equation}
\begin{equation}
    \mat{H} = 
    \sum_{m=1}^{M} \sum_{l=0}^{L-1}
    \big(\vu_l^{(m)}\big)^T \otimes
    \begin{bmatrix} \left(\frac{\vmuhat_{l+1}^{(m)}- \vmuhat_l^{(m)}}{\Delta t}\right)  & \left(\frac{\mSighat_{l+1,l}^{(m)} - \mSighat_{l,l}^{(m)}}{\Delta t} + \frac{\vmuhat_{l+1}^{(m)}- \vmuhat_l^{(m)}}{\Delta t}(\vmuhat_l^{(m)})^T\right) \end{bmatrix}.
\end{equation}
Assuming that $\mat{G}$ is positive-definite, let
\begin{equation}
    \begin{bmatrix} \mat{\tilde{V}}_0^T & \cdots & \mat{\tilde{V}}_{\dim\vect{u}}^T \end{bmatrix} =
    \mat{H}\mat{G}^{-1}.
\label{eqn:maximization_for_generators}
\end{equation}
With $\mat{A}_l^{(m)}$ and $\vect{b}_l^{(m)}$ defined by \cref{eqn:finite_time_Koopman_operator_approx} and 
$\mat{\tilde{C}} = \begin{bmatrix} \mat{I}_{\dim\vect{y}} & \mat{0} \end{bmatrix}$, let
\begin{multline}
    \mSig_{\vw} 
    = \frac{1}{M L}\sum_{m=1}^M\sum_{l=0}^{L-1}\Big[
    \mSighat_{l+1,l+1}^{(m)} 
    - \mat{A}_l^{(m)}\mSighat_{l,l+1}^{(m)}
    - \mSighat_{l+1,l}^{(m)}(\mat{A}_l^{(m)})^T
    + \mat{A}_l^{(m)}\mSighat_{l,l}^{(m)}(\mat{A}_l^{(m)})^T \\
    + \left(\vmuhat_{l+1}^{(m)} - \mat{A}_l^{(m)}\vmuhat_l^{(m)} - \vect{b}_l^{(m)}\right)\left(\vmuhat_{l+1}^{(m)} - \mat{A}_l^{(m)}\vmuhat_l^{(m)} - \vect{b}_l^{(m)}\right)^T
    \Big]
    + \frac{\gamma_{\mat{G}} \Delta t^2}{M L} \sum_{i=0}^{\dim\vu} \mat{\tilde{V}}_i^T \mat{\tilde{V}}_i
    + \frac{\gamma_{\vw}}{M L} \mat{I},
    \label{eqn:process_noise_covariance}
\end{multline}
\begin{equation}
    \vmu_0 = \frac{1}{M}\sum_{m=1}^M \vmuhat_0^{(m)}, 
\end{equation}
\begin{equation}
    \mSig_0 = \frac{1}{M}\sum_{m=1}^M\left[ \mSighat_{0,0}^{(m)} + \left(\vmuhat_0^{(m)} - \vmu_0\right)\left(\vmuhat_0^{(m)} - \vmu_0\right)^T \right] 
    + \frac{\gamma_0}{M} \mat{I},
    \label{eqn:initial_condition_covariance}
\end{equation}
\begin{equation}
    \vect{c}_0 
    = \frac{1}{M(L+1)} \sum_{m=1}^M \sum_{l=0}^L\left( \vect{y}_l^{(m)} - \mat{\tilde{C}}\vmuhat_l^{(m)} \right),
\end{equation}
\begin{multline}
    \mSig_{\vv} = \frac{1}{M(L+1)} \sum_{m=1}^M \sum_{l=0}^L \left[ \mat{\tilde{C}}\mSighat_{l,l}^{(m)}\mat{\tilde{C}}^T
    + \left(\vy_l^{(m)} - \vect{c}_0 - \mat{\tilde{C}}\vmuhat_l^{(m)}\right)\left(\vy_l^{(m)} - \vect{c}_0 - \mat{\tilde{C}}\vmuhat_l^{(m)}\right)^T \right] \\
    + \frac{\gamma_{\vv}}{M(L+1)} \mat{I}.
\end{multline}
The following theorem states that the above quantities are the unique maximizer of the regularized ELBO, provided that certain postitive-definiteness conditions hold.
\begin{theorem}[Maximization step]
\label{thm:maximization_step}
If $\mat{G}$, $\mSig_{\vw}$, $\mSig_0$, and $\mSig_{\vv}$ are positive-definite then $\calP = \big( \mat{\tilde{V}}_0, \ldots, \mat{\tilde{V}}_{\dim\vect{u}},\ \vect{c}_0,\ \mSig_{\vect{w}},\ \mSig_{\vect{v}},\ \vmu_0, \mSig_0 \big)$ is the unique maximizer of 
$\hat{L}_{Q} - R$, where $\hat{L}_{Q}$ is
the evidence lower bound given by \cref{eqn:ELBO}.
\end{theorem}
\begin{proof}
The proof is involved, so we give it in Appendix~\ref{app:maximization_step}.
\end{proof}

Nonzero coefficients in the regularization function are needed when the Gram matrix $\mat{G}$ or any of the covariance matrices $\mSig_{0}, \mSig_{\vv}, \mSig_{\vw}$ fail to be positive-definite, i.e., when they are singular.
This can happen when there is insufficient data or insufficient noise in the data, particularly when the dimension of $\vect{z}$ becomes large.
For example, if $M$ and $L$ are too small, then $\mat{G}$ can fail to be positive-definite.
Singularity can also occur when the model \cref{eqn:main_HMM} can fit some of the data perfectly.
Examining their definitions, the matrices $\mSig_{0}$, $\mSig_{\vv}$, and $\mSig_{\vw}$ measure how accurately the parameters $\mat{\tilde{V}}_0, \ldots, \mat{\tilde{V}}_{\dim\vu}$, $\vmu_0$, and $\vect{c}_0$ can be fit to the inferred latent state variables.
Hence, when some of the data can be fit without error, the covariance matrices $\mSighat^{(m)}_{l,l}$, $\mSig_{0}$, $\mSig_{\vv}$, or $\mSig_{\vw}$ can shrink to singularity during iteration of the EM algorithm.
Other methods to avoid singularity include restricting the model to use diagonal or block-diagonal $\mSig_{0}$, $\mSig_{\vv}$, or $\mSig_{\vw}$, or augmenting the training data with copies perturbed by small amounts of noise.
In the numerical applications presented later in \cref{sec:numerical_applications}, we train on enough data to prevent these singularities for the model dimensions we consider.

The solution for the generators during the maximization step of the EM algorithm provided by Theorem~\ref{thm:maximization_step} (with $R=0$) can be viewed as a limiting case of the control-affine extended dynamic mode decomposition technique presented in \cite{Peitz2020data}.
To see this, we observe that \cref{eqn:maximization_for_generators} can also be written as
\begin{multline}
    \begin{bmatrix} \mat{\tilde{V}}_0 \\ \vdots \\ \mat{\tilde{V}}_{\dim\vect{u}} \end{bmatrix} =
    \left( \frac{1}{M L}\sum_{m=1}^M \sum_{l=0}^{L-1} \mathbb{E}\left[ \left(\vect{u}_l^{(m)} \otimes \begin{bmatrix} 1 \\ \vect{z}_l^{(m)} \end{bmatrix} \right) \left(\vect{u}_l^{(m)} \otimes \begin{bmatrix} 1 \\ \vect{z}_l^{(m)} \end{bmatrix} \right)^T \ \Bigg\vert \ \mat{Y} = \mat{Y}^{(m)} \right] \right)^{-1} \\
    \left( \frac{1}{M L}\sum_{m=1}^M \sum_{l=0}^{L-1} \mathbb{E}\left[ \left(\vect{u}_l^{(m)} \otimes \begin{bmatrix} 1 \\ \vect{z}_l^{(m)} \end{bmatrix} \right) \left( \frac{\vect{z}_{l+1}^{(m)} - \vect{z}_{l}^{(m)} }{\Delta t} \right)^T \ \Bigg\vert \ \mat{Y} = \mat{Y}^{(m)} \right] \right),
    \label{eqn:maximization_for_generators_v2}
\end{multline}
where the conditional expectation is computed using the model parameters of the previous iteration thanks to \cref{eqn:exact_inference}.
To sample from these distributions, let trajectory indices $m_1, \ldots, m_K$ be drawn independently and uniformly from $\{ 1, \ldots, M \}$, let time indices $l_1, \ldots, l_k$ be drawn independently and uniformly from $\{ 0, \ldots, L-1 \}$ and draw $(\vpsi_k, \vpsi_k^{\sharp}) = (\vect{\hat{z}}_{l_k}^{(m_k)}, \vect{\hat{z}}_{l_k + 1}^{(m_k)})$ from the posterior distribution $P_{(\vect{z}_{l_k}, \vect{z}_{l_k+1})\vert \mat{Y} = \mat{Y}^{(m_k)}}$ for each $k=1,\ldots, K$.
Constructing matrices
\begin{equation}
    \boldsymbol{\Psi}_{K} = \frac{1}{\sqrt{K}}\begin{bmatrix} \vect{u}_{l_1}^{(m_1)} \otimes \begin{pmatrix} 1 \\ \vpsi_1 \end{pmatrix} & \cdots & \vect{u}_{l_K}^{(m_K)} \otimes \begin{pmatrix} 1 \\ \vpsi_K \end{pmatrix} \end{bmatrix},
\end{equation}
\begin{equation}
    \boldsymbol{\dot{\Psi}}_{K} = \frac{1}{\sqrt{K}}\begin{bmatrix}  \left(\frac{ \vpsi_1^{\sharp} - \vpsi_1 }{\Delta t}\right) & \cdots &  \left(\frac{ \vpsi_K^{\sharp} - \vpsi_K }{\Delta t}\right) \end{bmatrix},
\end{equation}
then yields the same solution as \cref{eqn:maximization_for_generators_v2} in the limit
\begin{equation}
    \begin{bmatrix} \mat{\tilde{V}}_0 \\ \vdots \\ \mat{\tilde{V}}_{\dim\vect{u}} \end{bmatrix} = \lim_{K\to\infty} \left( \boldsymbol{\Psi}_{K}\boldsymbol{\Psi}_{K}^T \right)^{-1}\boldsymbol{\Psi}_{K} \boldsymbol{\dot{\Psi}}_{K}^T,
\end{equation}
which holds almost surely thanks to the strong law of large numbers \cite{Koralov2012theory}.
Each term in the sequence are the generators computed by the method in \cite{Peitz2020data} from the snapshot pairs 
$\{ ( \vpsi_k, \vpsi_k^{\sharp} ) \}_{k=1}^K$.

It can also be readily shown that process noise covariance estimate is given by
\begin{equation}
    \frac{1}{\Delta t^2} \mSig_{\vw} = \frac{1}{M L}\sum_{m=1}^M \sum_{l=0}^{L-1} \mathbb{E}\left[ \left( \vect{e}_l^{(m)} \right) \left( \vect{e}_l^{(m)} \right)^T \ \Big\vert \ \mat{Y} = \mat{Y}^{(m)} \right],
\end{equation}
where
\begin{equation}
    \vect{e}_{l}^{(m)} 
    = \left(\frac{\vect{z}_{l+1}^{(m)} - \vect{z}_{l}^{(m)} }{\Delta t}\right) 
    - \left( \sum_{k=1}^{\dim\vu} [\vu_l^{(m)}]_k \mat{\tilde{V}}_k \right)^T \begin{bmatrix} 1 \\ \vect{z}_{l}^{(m)} \end{bmatrix}
\end{equation}
can be viewed as the modeling error for the dynamics of the states in our HMM.
From this expression, it is clear that the generators given by \cref{eqn:maximization_for_generators_v2} minimize the mean square error
\begin{multline}
    \frac{1}{\Delta t^2} \Tr(\mSig_{\vw}) = \\
    \frac{1}{M L}\sum_{m=1}^M \sum_{l=0}^{L-1} \mathbb{E} \left[ \left\Vert \frac{1}{\Delta t}\left(\begin{bmatrix} 1 \\ \vect{z}_{l+1}^{(m)} \end{bmatrix} - \begin{bmatrix} 1 \\ \vect{z}_{l}^{(m)} \end{bmatrix}\right) - \left( \sum_{k=1}^{\dim\vu} [\vu_l^{(m)}]_k \mat{V}_k \right)^T \begin{bmatrix} 1 \\ \vect{z}_{l}^{(m)} \end{bmatrix} \right\Vert^2 \ \Bigg\vert \ \mat{Y} = \mat{Y}^{(m)}  \right],
\end{multline}
recalling that $\mat{V}_k = \begin{bmatrix} \vect{0} & \mat{\tilde{V}}_k \end{bmatrix}$.
Upon convergence of the EM iterations, this quantity can be viewed as a residual measuring how closely the inferred latent states $\vect{z}$ obey input-affine dynamics.
In terms of the sample matrices constructed in the previous paragraph, the process noise covariance is obtained almost surely in the limit
\begin{equation}
    \frac{1}{\Delta t^2} \mSig_{\vw} = \lim_{K \to \infty} \boldsymbol{\dot{\Psi}}_{K} \boldsymbol{\dot{\Psi}}_{K}^T - 
    \boldsymbol{\dot{\Psi}}_{K} \boldsymbol{\Psi}_{K}^T \left( \boldsymbol{\Psi}_{K} \boldsymbol{\Psi}_{K}^T \right)^{-1} \boldsymbol{\Psi}_{K} \boldsymbol{\dot{\Psi}}_{K}^T,
\end{equation}
with its trace given by the residual
\begin{equation}
    \frac{1}{\Delta t^2} \Tr(\mSig_{\vw}) 
    = \min_{\mat{\tilde{V}}_0, \ldots, \mat{\tilde{V}}_{\dim\vect{u}} } \lim_{K\to\infty} \left\Vert \boldsymbol{\dot{\Psi}}_{K} - \begin{bmatrix} \mat{\tilde{V}}_0^T & \cdots & \mat{\tilde{V}}_{\dim\vu}^T \end{bmatrix} \boldsymbol{\Psi}_{K} \right\Vert_F^2.
\end{equation}
This closely resembles the residual developed in \cite{Colbrook2021rigorous, Colbrook2023residual} to measure the spectral leakage when the Koopman operator is projected onto finite-dimensional spaces of observables.
Making this connection rigorous will be an important topic of future work.

Another important observation is that the explicit solutions for the parameters during the maximization step of the EM algorithm given in Theorem~\ref{thm:maximization_step} depend only on a few parameters of the inference distribution.
Therefore, during the expectation or E-step of the EM algorithm, we need only compute the conditional or ``posterior'' mean and covariance
\begin{equation}
    \vmuhat_l^{(m)} = \mathbb{E}\left[ \vect{z}_l\ \vert \ \mat{Y} = \mat{Y}^{(m)} \right], \qquad
    \mSighat_{l,l}^{(m)} = \mathbb{E}\left[ \left(\vect{z}_l - \vmuhat_l^{(m)}\right)\left(\vect{z}_{l} - \vmuhat_{l}^{(m)}\right)^T \ \vert \ \mat{Y} = \mat{Y}^{(m)} \right]
    \label{eqn:posterior_means_variances}
\end{equation}
at each time step given the observations $\mat{Y}^{(m)}$ along each trajectory as well as the posterior covariance between adjacent time steps
\begin{equation}
    \mSighat_{l,l+1}^{(m)} = \mathbb{E}\left[ \left(\vect{z}_l - \vmuhat_l^{(m)}\right)\left(\vect{z}_{l+1} - \vmuhat_{l+1}^{(m)}\right)^T \ \vert \ \mat{Y} = \mat{Y}^{(m)} \right].
    \label{eqn:eqn:adjacent_posterior_variances}
\end{equation}
Fortunately, there is a very efficient algorithm due to R. H. Shumway and D. S. Stoffer \cite{Shumway1982approach} for computing these conditional expectations that proceeds by passing over each trajectory twice.
On the first pass, we move forward along the trajectory using a Kalman filter \cite{Kalman1960new} to assimilate the observations made up to each given time step.
On the second pass, we move backward along the trajectory using a smoother \cite{Rauch1963solutions} to assimilate the observations made after each given time step.
The details of this ``forward-backward algorithm'' are provided in Section~\ref{subsec:Estep}.


\subsection{Expectation step: Kalman filter and smoother}
\label{subsec:Estep}
In the expectation step of the EM algorithm, we use the data to estimate the latent states $\vect{z}_l$ along each trajectory, given fixed model parameters~$\mathcal{P}$.  Recall that with the parameters $\mathcal{P}$ fixed, the dynamics of our discrete-time model \cref{eqn:main_HMM} can be written as
\begin{equation}
\begin{split}
    \vect{z}_{l+1} &= \mat{A}_l\vect{z}_l + \vect{b}_l + \vect{w}_l \\
    \vect{y}_l &= \vect{c}_0 + \mat{\tilde{C}}\vect{z}_l + \vect{v}_l,
\end{split}
\label{eqn:E_step_dynamical_model}
\end{equation}
where $\mat{A}_l$, $\vect{b}_l$, $\vect{c}_0$, and $\mat{\tilde{C}}$ are all known. 
Furthermore, the initial condition, process noise, and measurement noise have independent Gaussian distributions
\begin{equation}
    \vect{z}_0\sim \mathcal{N}(\vmu_0, \mSig_0), \qquad
    \vect{w}_l \sim \mathcal{N}(\vect{0}, \mSig_{\vect{w}}), \qquad
    \vect{v}_l \sim \mathcal{N}(\vect{0}, \mSig_{\vect{v}}),
\end{equation}
with known means and covariances.
Therefore, computing the posterior means and covariances \cref{eqn:posterior_means_variances}~and~\cref{eqn:eqn:adjacent_posterior_variances} of the states $\vect{z}_l$ given the observations $\mathbf{Y}^{(m)}$ along the $m$th independent trajectory is a standard state estimation problem for the linear time-varying dynamical system \cref{eqn:E_step_dynamical_model} that can be solved by means of Kalman filtering and smoothing \cite{Kalman1960new, Rauch1963solutions, Anderson1979optimal, Rhodes1971tutorial}.

This section summarizes the appropriate filtering and smoothing equations, which are essentially the same as those given in \cite{Yu2004derivation, Ghahramani1996parameter, Shumway1982approach, Rauch1963solutions, Anderson1979optimal}, but not identical.
We derive these equations in \cref{app:expectation_step} by closely following the derivation in \cite{Yu2004derivation}.

During the forward pass, the posterior covariance of each state $\vect{z}_l$ given the observations
$\vect{y}_0, \ldots, \vect{y}_l$ up to the $l$th time step is computed using an efficient recursion process.
Let us denote posterior means and covariances by
\begin{equation}
    \vmuhat_{k\vert l} = \mathbb{E}\left[ \vect{z}_k\ \vert \ \vect{y}_0, \ldots, \vect{y}_l \right], \qquad
    \mSighat_{j,k\vert l} = \mathbb{E}\left[ (\vect{z}_j - \vmuhat_{j\vert l} ) (\vect{z}_k - \vmuhat_{k\vert l} )^T\ \vert \ \vect{y}_0, \ldots, \vect{y}_l \right]
\end{equation}
and keep in mind that we are ultimately looking for $\vmuhat_l = \vmuhat_{l\vert L}$, $\mSighat_{l,l} = \mSighat_{l,l\vert L}$, and $\mSighat_{l,l+1} = \mSighat_{l,l+1\vert L}$.
Initially, given no observations from the trajectory, we set
\begin{equation}
    \vmuhat_{0\vert -1} = \vmu_0, \qquad \mSighat_{0,0\vert -1} = \mSig_0.
\end{equation}
On the other hand, if we have the previous conditional covariance $\mSighat_{l-1,l-1\vert l-1}$, then the next covariance given the same set of observations is
\begin{equation}
    \mSighat_{l,l\vert l-1} = \mSig_{\vect{w}} + \mat{A}_{l-1}\mSighat_{l-1,l-1\vert l-1}\mat{A}_{l-1}^T.
    \label{eqn:forward_pass_covariance_increment}
\end{equation}
These covariances are then used to update the conditional covariance at the current step given the observations so far using
\begin{equation}
    \mSighat_{l,l\vert l} = \mSighat_{l,l\vert l-1} - \mSighat_{l,l\vert l-1}\mat{\tilde{C}}^T\left( \mSig_{\vect{v}} + \mat{\tilde{C}} \mSighat_{l,l\vert l-1} \mat{\tilde{C}}^T \right)^{-1}\mat{\tilde{C}}\mSighat_{l,l\vert l-1}.
    \label{eqn:forward_pass_covariance_update}
\end{equation}
The process is repeated in order to estimate each $\mSighat_{l,l\vert l}$ recursively.
Given the previous estimate of the mean $\vmuhat_{l-1\vert l-1}$, the so called ``Kalman gain''
\begin{equation}
    \mat{K}_{l} = \mSighat_{l,l\vert l-1} \mat{\tilde{C}}^T \left( \mSig_{\vect{v}} + \mat{\tilde{C}} \mSighat_{l,l\vert l-1} \mat{\tilde{C}}^T \right)^{-1}
    \label{eqn:Kalman_gain}
\end{equation}
is used to assimilate the latest measurement into our estimate of the current mean via
\begin{equation}
    \vmuhat_{l\vert l} = \vmuhat_{l\vert l-1} + \mat{K}_{l}\left( \vect{y}_l - \vect{c}_0 - \mat{\tilde{C}}\vmuhat_{l\vert l-1} \right), \qquad
    \vmuhat_{l\vert l-1} = \mat{A}_{l-1}\vmuhat_{l-1\vert l-1} + \vect{b}_{l-1}.
    \label{eqn:forward_pass_mean_update}
\end{equation}
During the forward pass, we store each $\vmuhat_{l\vert l}$, $\mSighat_{l,l\vert l}$, and $\mSighat_{l,l\vert l-1}$ to be used during the smoothing step.

The backward pass or ``smoother'' incorporates the remaining observations $\vy_{l+1}, \ldots, \vy_L$ into the optimal estimate of each $\vz_l$ using another recursion.
Given the mean of the posterior Gaussian distribution one time-step in the future $\vmuhat_{l+1}$, the smoothing gain
\begin{equation}
    \mat{J}_l = \mSighat_{l,l\vert l} \mat{A}_l^T\left( \mSighat_{l+1,l+1\vert l} \right)^{-1}
    \label{eqn:Kalman_smoothing_gain}
\end{equation}
is computed in order to update the mean of the posterior distribution at the current time-step
\begin{equation}
    \vmuhat_l = \vmuhat_{l\vert l} + \mat{J}_l\left(\vmuhat_{l+1} - \mat{A}_l \vmuhat_{l\vert l} - \vect{b}_l \right).
    \label{eqn:backward_pass_mean_update}
\end{equation}
The posterior covariance matrices are also computed recursively according to
\begin{equation}
    \mSighat_{l,l} = \mSighat_{l,l\vert l} + \mat{J}_l\left(\mSighat_{l+1,l+1} - \mSighat_{l+1,l+1\vert l}  \right) \mat{J}_l^T, \qquad
    \mSighat_{l,l+1} = \mat{J}_l\mSighat_{l+1,l+1}.
    \label{eqn:backward_pass_covariance_update}
\end{equation}
The backward recursion is initialized with the final step of the Kalman filter $\vmuhat_L = \vmuhat_{L\vert L}$, and $\mSighat_{L,L} = \mSighat_{L,L\vert L}$.
Finally, we have obtained the required parameters of the inference distribution, namely $\lbrace \vmuhat_l \rbrace_{l=0}^L$, $\lbrace \mSighat_{l,l} \rbrace_{l=0}^L$, and $\lbrace \mSighat_{l,l+1} \rbrace_{l=0}^{L-1}$, to be used during the maximization stage of the EM algorithm.

\section{Approximating Koopman eigenfunctions}
\label{sec:approximating_eigenfunction}

Once we have learned the parameters $\mathcal{P}$, we use the model \cref{eqn:finite_dimensional_generator_IO_model} as a surrogate for the original dynamical system.
As discussed in \cite{Mezic-05, Mezic:2019, Mauroy2016global, Mauroy2013isostables}, obtaining accurate approximations of Koopman eignenfunctions reveals useful information about the system.
It is easy to compute the values of eigenfunctions of the Koopman generators associated with the surrogate model when the state $\vect{z} = \vpsi(\vect{x})$ and input $\vect{u}$ are known.
Specifically, if $\vect{v}$ is an eigenvector with eigenvalue $\lambda$ of the matrix 
\begin{equation}
    \mat{V}_{\vect{u}} = \sum_{k=1}^{\dim\vect{u}} u_k \mat{V}_k,
\end{equation}
then the function $\varphi = \vect{v}^T\vpsi$ evolves under the surrogate model \cref{eqn:finite_dimensional_generator_IO_model} according to
\begin{equation}
    \frac{d}{d t} \varphi(\vect{x}) 
    = \vect{v}^T \frac{d}{d t}\vpsi(\vect{x}) 
    = \vect{v}^T \mat{V}_{\vect{u}}^T \vpsi(\vect{x})
    = \lambda \vect{v}^T\vpsi = \lambda \varphi(\vect{x}).
\end{equation}
That is, $\varphi$ is an eigenfunction with eigenvalue $\lambda$ of the Koopman generator for the surrogate model with input $\vu$.

Since we do not have direct access to the values $\vect{z} = \vpsi(\vect{x})$, we use the conditional expectation of $\vect{z}$ given the observed outputs $\vect{y}_0, \ldots, \vect{y}_L$ along a trajectory.
The filtering and smoothing algorithm presented in \cref{subsec:Estep} allows us to compute
\begin{equation}
    \vmuhat_l = \mathbb{E}\left[ \vect{z}_l \ \vert \ \mat{Y} = \begin{bmatrix} \vect{y}_0, \ldots, \vect{y}_L \end{bmatrix} \right].
\end{equation}
We use this as an estimate of $\vpsi(\vect{x}(t_l))$ when computing the values of eigenfunctions of the surrogate model.
That is, we compute $\vect{v}^T \vmuhat_l$ as an estimate of $\varphi(\vect{x}(t_l))$.

\section{Numerical applications}
\label{sec:numerical_applications}
In this section, we study the performance of our proposed method in different areas, namely, the identification of Koopman eigenfunctions, the prediction quality for observed quantities $\vect{y}$, and the control performance when using \cref{eqn:finite_dimensional_generator_IO_model} as a surrogate model for feedback control in a model predictive control (MPC) framework.

As a control objective, we aim to track a reference trajectory for the output~$\vect{y}$. 
Within MPC, we need to solve an open-loop optimal control problem over a prediction horizon of finite length $n_p \Delta t$. We then apply a fraction of the optimal control signal $u^*$ to the actual system over the control horizon of length $n_c \Delta t$ with $n_c \leq n_p$. 
To realize this in real-time, we need to solve an optimization problem of the form
\begin{equation}\label{eqn:OCP_original}
\begin{aligned}
\minimize_{u} \int_t^{t+n_p\Delta t} \Big[\big(\vy-\vy^{\mathsf{ref}}\big)^T &\mQ \big(\vy-\vy^{\mathsf{ref}}\big) + \vu^T \mR \vu\Big] \, dt \\
\mbox{s.t.} \quad \frac{d}{dt}\vx &= \sum_{i=0}^q u_i \vf_i(\vx), \quad \vx(t) = \vx_0 \\
\vy &= \vg(\vx) + \vect{w}_t, \\
\vu^{\mathsf{min}} &\leq \vu \leq \vu^{\mathsf{max}},
\end{aligned}
\end{equation}
at each step of MPC.
Here, $\mQ$ and $\mR$ are positive semidefinite and positive definite matrices of appropriate sizes, respectively, and $\vx_0$ is the initial condition that is obtained by measuring the state of the plant. 
As problem~\cref{eqn:OCP_original} may be too costly for real-time control if
the original dynamics \cref{eq:bilinear_system} are complex (or may even be unknown),
we use the surrogate model \cref{eqn:finite_dimensional_generator_IO_model} instead of the original differential equation, which results in the following problem:
\begin{equation}\label{eqn:OCP_surrogate}
\begin{aligned}
\minimize_{u} \Delta t \sum_{j=k}^{k+n_p} \Big[\big(\vy_j-\vy_j^{\mathsf{ref}}\big)^T &\mQ \big(\vy_j-\vy_j^{\mathsf{ref}}\big) + \vu_j^T \mR \vu_j \Big] \\
\mbox{s.t.} \qquad\qquad \mbox{\cref{eqn:finite_dimensional_generator_IO_model} is satisfied} \quad \mbox{and}\quad
&\vu^{\mathsf{min}} \leq \vu \leq \vu^{\mathsf{max}}.
\end{aligned}
\end{equation}
To estimate the initial condition for the latent state $\vz_k$ from the history of incoming measurements $\{\vy_0, \ldots, \vy_k\}$ and past control inputs $\{\vu_0, \ldots, \vu_{k-1}\}$, the Kalman filter as described in Section \ref{subsec:Estep} is used in each loop of the MPC routine. 
Problem \cref{eqn:OCP_surrogate} can now be solved using efficient gradient-based optimization, where we calculate the derivative of the objective function via the adjoint equation corresponding to \cref{eqn:finite_dimensional_generator_IO_model} (see \cite{Peitz2020data} for a detailed derivation).

The EM algorithm requires an initial guess for the matrices $\mat{V}_i$.  If there is no control input (as in the example in \cref{sec:duffing} below), then we use EDMD to provide the initial guess.  If a control input is present, as in the other examples below, then we initialize with random matrices chosen 
so that the eigenvalues of $\mat{I} + \Delta t \mat{V}_0$ and each matrix $\Delta t \max\{\vert u_i^{\mathsf{min}} \vert, \vert u_i^{\mathsf{max}} \vert \} \mat{V}_i$, $i \geq 1$ are approximately distributed across the unit disk according to Girko's circular law \cite{Girko1985circular}.
Because the results depend on the initialization, we use multiple random initializations and we take the model with the best performance.

\subsection{Actuated system with a polynomial slow manifold}
We consider the following system from \cite{Otto2021koopman} that was originally adapted from an example in \cite{Brunton2016koopman}:
\begin{equation}
\begin{split}
\dot{x}_1 &= -\alpha x_1 + u \\
\dot{x}_2 &= \beta\left( x_1^3 - x_2 \right) \\
y &= x_2 + v.
\end{split}
\end{equation}
This system can be fully described in a finite-dimensional Koopman-invariant subspace as
\begin{equation}
\begin{split}
\frac{d}{dt}\begin{bmatrix} 1 \\ x_1 \\ x_2 \\ x_1^2 \\ x_1^3 \end{bmatrix} &= 
\begin{bmatrix} 
0 & 0 & 0 & 0 & 0 \\
0 & -\alpha & 0 & 0 & 0 \\
0 & 0 & -\beta & 0 & \beta \\
0 & 0 & 0 & -2\alpha & 0 \\
0 & 0 & 0 & 0 & -3\alpha 
\end{bmatrix}
\begin{bmatrix} 1 \\ x_1 \\ x_2 \\ x_1^2 \\ x_1^3 \end{bmatrix}
+ u \begin{bmatrix} 
0 & 0 & 0 & 0 & 0 \\
1 & 0 & 0 & 0 & 0 \\
0 & 0 & 0 & 0 & 0 \\
0 & 2 & 0 & 0 & 0 \\
0 & 0 & 0 & 3 & 0
\end{bmatrix}
\begin{bmatrix} 1 \\ x_1 \\ x_2 \\ x_1^2 \\ x_1^3 \end{bmatrix} \\
y &= x_2 + v.
\end{split}
\end{equation}
We choose $\alpha = 1$ and $\beta = 5$ with Gaussian measurement noise $v$ with zero mean and variance $0.01$, and zero process noise.
The inputs are held constant over time intervals of length $0.5$ and take values drawn from a Gaussian distribution with zero mean and variance $\sigma_u^2=5$.

We collected $50$ independent trajectories with $500$ observations recorded at intervals $\Delta t = 0.01$.
The initial conditions were drawn from an isotropic Gaussian distribution with zero mean and unit variance.
The first $250$ points along each trajectory were used for training while the rest were saved for testing the model.

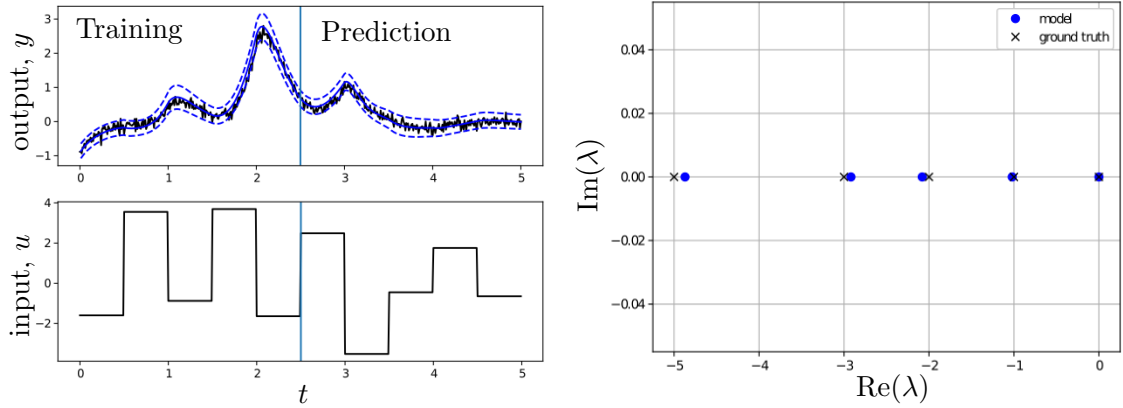
\begin{figure}
	\centering
	\begin{tikzonimage}[trim=0 0 0 0, clip=true, width=0.455\textwidth]{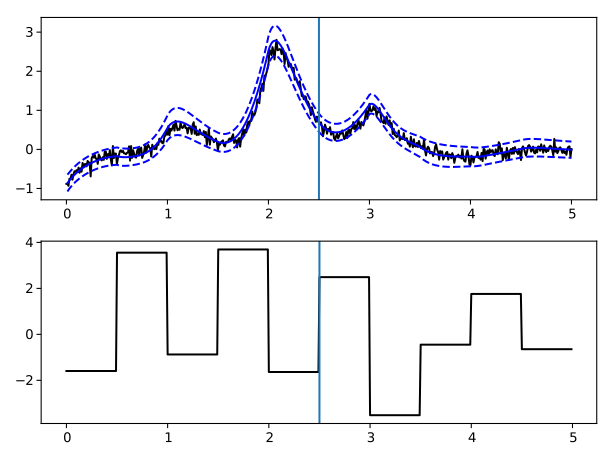}
		\node[rotate=0] at (0.525, 0.00) {$t$};
		\node[rotate=90] at (0.0, 0.75) {output, $y$};
		\node[rotate=90] at (0.0, 0.27) {input, $u$};
		\node at (0.2,0.9) {Training};
		\node at (0.68,0.9) {Prediction};
	\end{tikzonimage}
	\begin{tikzonimage}[trim=10 10 10 10, clip=true, width=0.49\textwidth]{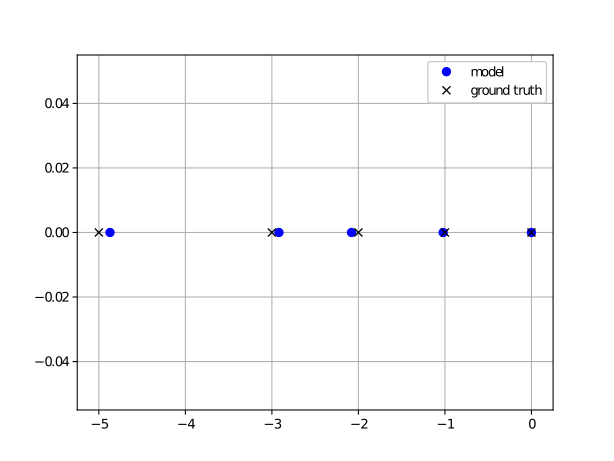}
		\node[rotate=0] at (0.525, 0.00) {Re($\lambda$)};
		\node[rotate=90] at (0.0, 0.5) {Im($\lambda$)};
	\end{tikzonimage}
	\caption{Left: Observations from the actuated toy model (black lines) together with our model's prediction (blue line) and $2\sigma$ confidence interval (dashed blue lines). The initial condition for the model prediction was computed using optimal state estimation from the training interval (left of vertical line). Right: Eigenvalues of the matrix approximation for the drift Koopman generator learned by the EM algorithm for the toy model.}
	\label{fig:toy_model_preds}
\end{figure}


The correct model dimension $\dim(\vect{z}) = 5$ is revealed by the largest likelihood among converged models with different dimensions and several random initializations.
After training the model, we used it to predict a trajectory and its uncertainty shown in \Cref{fig:toy_model_preds} on the left.
This was done by using the training portion of the trajectory to obtain an optimal estimate of the initial condition as well as its uncertainty.
The mean and uncertainty were then propagated forward by the model dynamics from the initial condition in order to produce the predicted observations and uncertainty shown in \Cref{fig:toy_model_preds}.
The eigenvalues of the matrix approximation of the zero input or ``drift'' Koopman generator learned by our model are compared to the ground truth eigenvalues in \cref{fig:toy_model_preds} on the right, showing excellent agreement.

\subsection{Duffing oscillator}
\label{sec:duffing}
The next system we study is the (unforced) Duffing oscillator:
\begin{equation}
\label{eq:Duffing}
\Ddot{x} + \delta \dot{x} + x(\beta + \alpha x^2) = 0, \qquad
\vect{y} = \vect{x} = (x, \dot{x})
\end{equation} 
with $\alpha = -1$, $\beta = 1$ and $\delta = 0.5$, which are the parameters that were also studied in \cite{Williams2015data,Otto2019linearly}. 
With these parameters, the system has stable fixed points at $x=\pm 1$ with eigenvalues $\lambda_{1,2} = \frac{1}{4}\big( -1 \pm i \sqrt{31} \big)$ and a saddle point at the origin with eigenvalues $\lambda_{3,4} = \frac{1}{4}\big( -1 \pm \sqrt{17} \big)$.
The eigenvalues of the linearizations about the stable spiral fixed points are also eigenvalues of the Koopman generator.
The magnitude and phase of the corresponding complex eigenfunctions provide polar coordinates in each basin of attraction \cite{Lan2013linearization,Williams2015data}.
There are also Koopman eigenfunctions that take different constant values in the two basins of attraction of the stable fixed points.

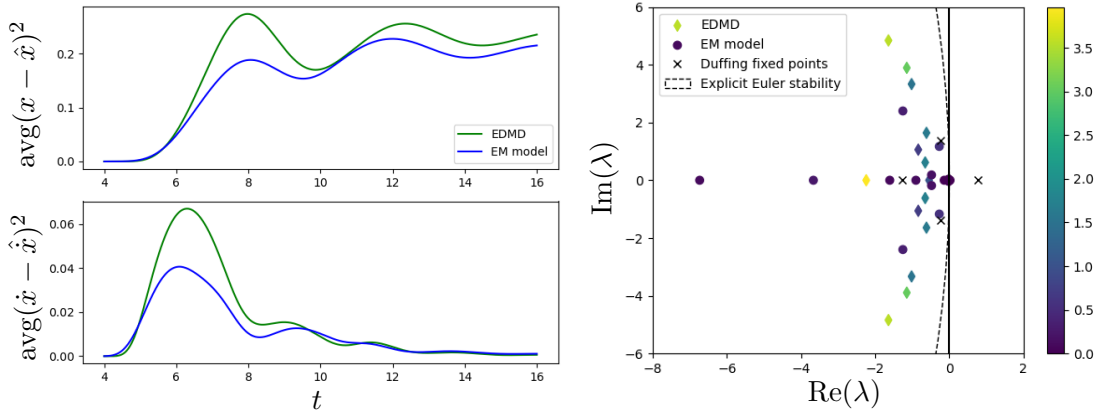
\begin{figure}
	\centering
	\begin{tikzonimage}[trim=0 0 0 0, clip=true, width=0.454\textwidth]{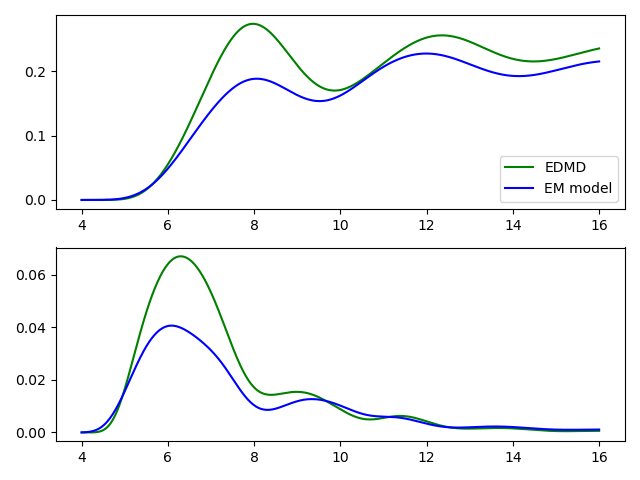}
		\node[rotate=0] at (0.525, -0.01) {$t$};
		\node[rotate=90] at (-0.02, 0.75) {$\avg (x - \hat{x})^2$};
		\node[rotate=90] at (-0.02, 0.27) {$\avg (\dot{x} - \hat{\dot{x}})^2$};
	\end{tikzonimage}
	\begin{tikzonimage}[trim=20 10 10 10, clip=true, width=0.475\textwidth]{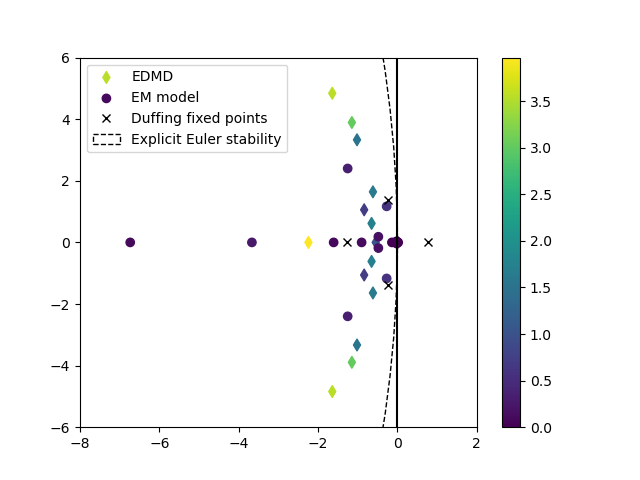}
		\node[rotate=0] at (0.43, -0.01) {Re($\lambda$)};
		\node[rotate=90] at (0.0, 0.5) {Im($\lambda$)};
	\end{tikzonimage}
	\caption{Left: Mean square prediction error for the Duffing equation using the model found by EDMD and the model trained using the proposed EM algorithm. The initial condition for the EM model prediction was computed using optimal state estimation over the time interval $0\leq t < 4$. Right: Eigenvalues of the matrix approximations for the Koopman generator found by EDMD and the EM algorithm for the Duffing equation. The eigenvalues are colored by the residuals $\res_{M,\Delta t}(\lambda, \varphi)$, where $\varphi$ are the approximate eigenfunctions.}
	\label{fig:Duffing_error_and_eigenvalues}
\end{figure}

For simplicity we observe the full state $\vect{y} = \vect{x} = (x, \dot{x})$ without noise.
For training, we collected $801$ equally spaced samples along each of 50 trajectories with $0 \leq t \leq 16$, $\Delta t = 0.02$, and initial conditions $(x(0), \dot{x}(0))$ drawn uniformly at random from the box $[-2,2]^2$.
The testing data were collected in the same way.
We initialized the EM algorithm using the $16$-dimensional model found by EDMD using a dictionary of multivariate Legendre polynomials orthogonal on the box $[-2,2]^2$ with degrees up to $3$ in $x$ and $\dot{x}$.
The EM algorithm converged in the sense that it failed to increase the likelihood after $55$ iterations.
We see from the left plot in \cref{fig:Duffing_error_and_eigenvalues} that the EM algorithm reduced the mean square prediction error along trajectories when compared to EDMD.
Here, we estimated the HMM's state at time $t=4$ based on the observed testing data over the interval $0\leq t < 4$, and subsequent predictions are based only on these observations.
Unlike EDMD, the hidden Markov model's state is not an explicit function of the system's state and must be estimated based on an observed time history.
Hence, the prediction accuracy depends on the length of the interval used to estimate the model's state, with longer intervals resulting in more accurate state estimates and forecasts.

We use the residual proposed in \cite{Colbrook2021rigorous, Colbrook2023residual} in order to evaluate the quality of the approximate eigenvalue-eigenfunction pairs $(\lambda, \varphi)$ computed using EDMD and the EM algorithm.
The approximate eigenfunction values $\varphi(\vect{x}_l^{(m)}) = \vect{v}^T \vmuhat_l^{(m)}$ for the HMM at the $l$th sample along the $m$th training trajectory were obtained using the method described in \cref{sec:approximating_eigenfunction}.
Following \cite{Colbrook2021rigorous, Colbrook2023residual}, we define an empirical inner product
\begin{equation}
    \langle f, g \rangle_{M} = \frac{1}{M L} \sum_{m=1}^M \sum_{l=0}^{L-1} \bar{f}\big(\vect{x}_l^{(m)}\big) g\big(\vect{x}_l^{(m)}\big),
\end{equation}
and the associated semi-norm $\Vert f \Vert_{M} = \sqrt{\langle f, f \rangle_{M}}$ using the the training data.
The empirical residual is then defined by
\begin{equation}
    \res_{M, \Delta t}(\lambda, \varphi) 
    = \frac{\big\Vert \frac{1}{\Delta t}(U^{\Delta t} - I) \varphi - \lambda \varphi \big\Vert_{M}}{\Vert \varphi \Vert_{M}},
\end{equation}
where $U^{\Delta t}$ is the time-$\Delta t$ Koopman operator of \cref{eq:Duffing}.
At each sample point, the action of the Koopman operator on a function is given by $U^{\Delta t} \varphi( \vect{x}_l^{(m)} ) = \varphi( \vect{x}_{l+1}^{(m)} )$.
This allows us to compute the above residual from the collected data using the method described in \cite{Colbrook2023residual}.
The limiting residual $\lim_{\Delta t \to 0} \lim_{M\to\infty} \res_{M, \Delta t}(\lambda, \varphi)$ measures whether the function $\varphi$ is truly an eigenfunction of the Koopman generator $V$ with eigenvalue $\lambda$.


\begin{table}[]
    \centering
    \begin{tabular}{|c|c|c|c|c|}
        \hline
        \multirow{2}{4em}{EDMD} & $\lambda$ & $-0.002500$ & $-0.8387\pm 1.059 i$ & $-1.019 \pm 3.331 i$ \\
        & $\res_{M,\Delta t}(\lambda, \varphi)$ & $0.2722$ & $0.7061$ & $1.528$ \\
        \hline
        \multirow{2}{4em}{EM} & $\lambda$ & $ 0.01871$ & $-0.2720 \pm 1.170 i$ & $-1.254 \pm 2.400 i$ \\
        & $\res_{M,\Delta t}(\lambda, \varphi)$ & $0.008380$ & $0.5629$ & $0.3391$ \\
        \hline
    \end{tabular}
    \vspace{1em}
    \caption{Selected eigenvalues and residuals obtained using EDMD and the EM algorithm for the Duffing oscillator.}
    \label{tab:Duffing_EDMD_eigs_and_resids}
\end{table}

The right plot in \cref{fig:Duffing_error_and_eigenvalues} shows the eigenvalues of the Koopman generator approximations obtained using EDMD and the EM algorithm compared to the eigenvalues of the Duffing equation linearized about its fixed points.
The approximate eigenvalues are colored according to their empirical residuals.
The eigenvalues of the HMM have smaller residuals than all of the eigenvalues computed by EDMD with the exception of the EDMD eigenvalue $-0.002500$.
The learned HMM has an eigenvalue close to zero, as expected in order to described the two invariant basins of the Duffing equation.
It also has complex conjugate pairs of eigenvalues close to the expected eigenvalues $\lambda_{1,2}$ and their harmonics $2 \lambda_{1,2}$.
We tabulate these eigenvalues, their residuals, and nearby EDMD eigenvalues with low residuals in \cref{tab:Duffing_EDMD_eigs_and_resids}.

\begin{figure}
	\centering
	\begin{tikzonimage}[trim=0 0 0 40, clip=true, width=0.45\textwidth]{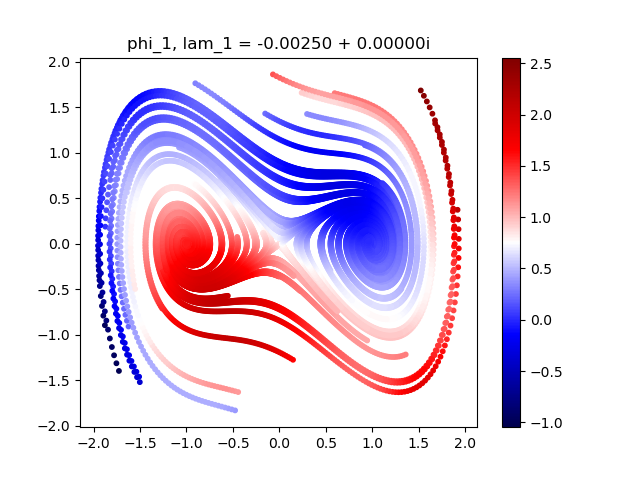}
	    \node[rotate=0] at (0.44, 1.075) {\footnotesize $\varphi$, $\lambda=0.00250$};
	    \node[rotate=90] at (0.02, 0.55) {\footnotesize $\dot{x}$};
	    \node[rotate=0] at (0.44, 0.02) {\footnotesize $x$};
	\end{tikzonimage}
	\begin{tikzonimage}[trim=0 0 0 40, clip=true, width=0.45\textwidth]{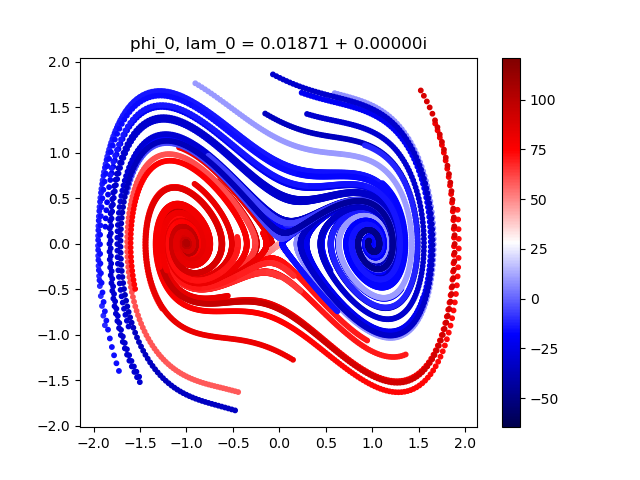}
	    \node[rotate=0] at (0.44, 1.075) {\footnotesize $\varphi$, $\lambda=0.01871$};
	    \node[rotate=90] at (0.02, 0.55) {\footnotesize $\dot{x}$};
	    \node[rotate=0] at (0.44, 0.02) {\footnotesize $x$};
	\end{tikzonimage}
	\caption{Training data for the Duffing equation colored by the approximate values of the eigenfunction with eigenvalue closest to zero, obtained by EDMD (left), and the EM algorithm (right). }
	\label{fig:Duffing_eigenfunction0}
\end{figure}

\begin{figure}
	\centering
	\begin{tikzonimage}[trim=0 0 0 40, clip=true, width=0.45\textwidth]{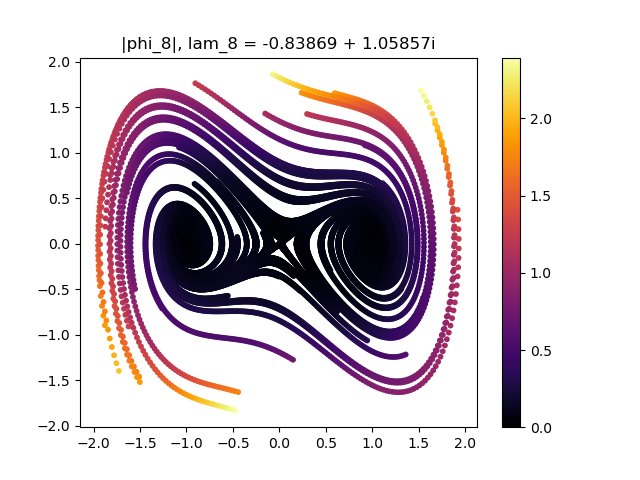}
	    \node[rotate=0] at (0.44, 1.075) {\footnotesize $\vert\varphi\vert$, $\lambda = -0.8387\pm 1.059 i$};
	    \node[rotate=90] at (0.02, 0.55) {\footnotesize $\dot{x}$};
	    \node[rotate=0] at (0.44, 0.02) {\footnotesize $x$};
	\end{tikzonimage}
	\begin{tikzonimage}[trim=0 0 0 40, clip=true, width=0.45\textwidth]{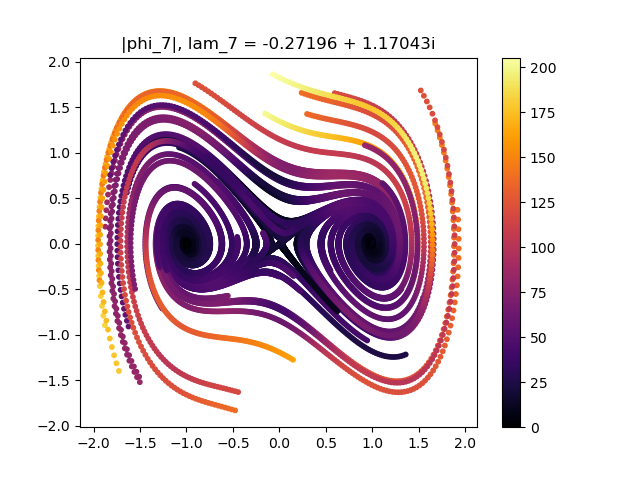}
	    \node[rotate=0] at (0.44, 1.075) {\footnotesize $\vert\varphi\vert$, $\lambda = -0.2720 + 1.170 i$};
	    \node[rotate=90] at (0.02, 0.55) {\footnotesize $\dot{x}$};
	    \node[rotate=0] at (0.44, 0.02) {\footnotesize $x$};
	\end{tikzonimage}\\
	\begin{tikzonimage}[trim=0 0 0 40, clip=true, width=0.45\textwidth]{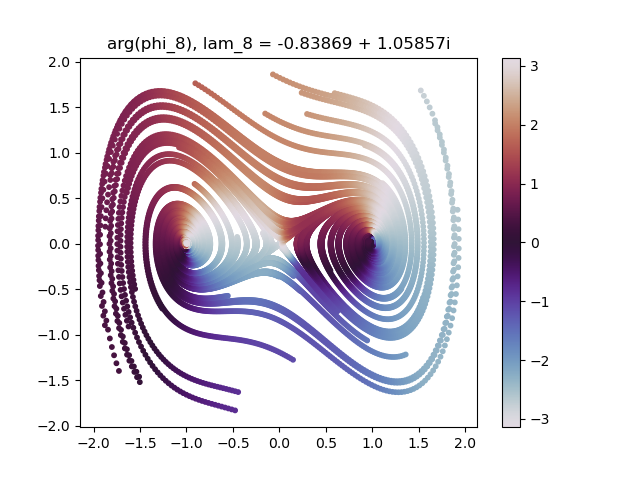}
	    \node[rotate=0] at (0.44, 1.075) {\footnotesize $\angle \varphi$, $\lambda = -0.8387\pm 1.059 i$};
	    \node[rotate=90] at (0.02, 0.55) {\footnotesize $\dot{x}$};
	    \node[rotate=0] at (0.44, 0.02) {\footnotesize $x$};
	\end{tikzonimage}
	\begin{tikzonimage}[trim=0 0 0 40, clip=true, width=0.45\textwidth]{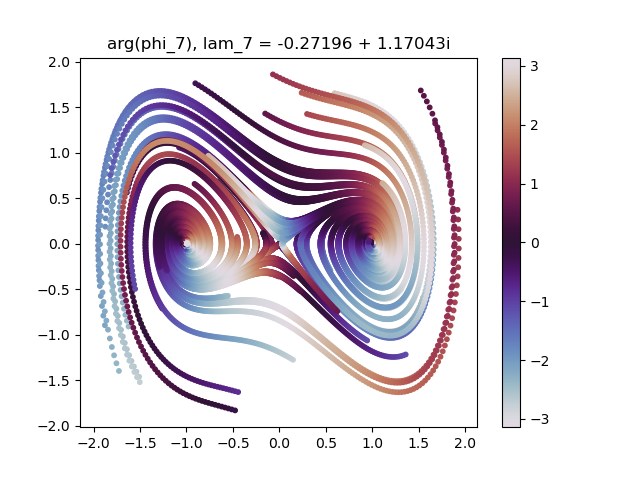}
	    \node[rotate=0] at (0.44, 1.075) {\footnotesize $\angle \varphi$, $\lambda = -0.2720 + 1.170 i$};
	    \node[rotate=90] at (0.02, 0.55) {\footnotesize $\dot{x}$};
	    \node[rotate=0] at (0.44, 0.02) {\footnotesize $x$};
	\end{tikzonimage}
	\caption{Analog of \cref{fig:Duffing_eigenfunction0} for the eigenvalue closest to that of the linearization about the stable fixed points: EDMD (left); EM algorithm (right). }
	\label{fig:Duffing_eigenfunction1}
\end{figure}

\begin{figure}
	\centering
	\begin{tikzonimage}[trim=0 0 0 40, clip=true, width=0.45\textwidth]{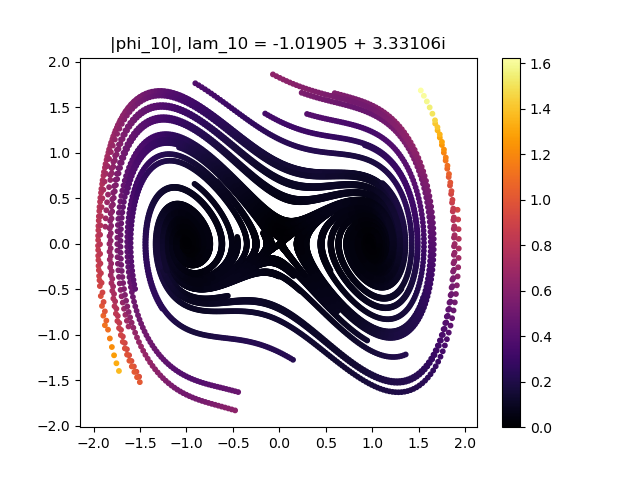}
	    \node[rotate=0] at (0.44, 1.075) {\footnotesize $\vert\varphi\vert$, $\lambda = -1.019 + 3.331 i$};
	    \node[rotate=90] at (0.02, 0.55) {\footnotesize $\dot{x}$};
	    \node[rotate=0] at (0.44, 0.02) {\footnotesize $x$};
	\end{tikzonimage}
	\begin{tikzonimage}[trim=0 0 0 40, clip=true, width=0.45\textwidth]{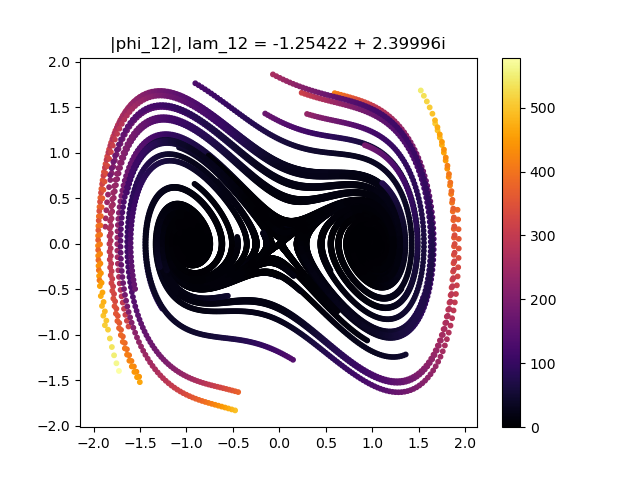}
	    \node[rotate=0] at (0.44, 1.075) {\footnotesize $\vert\varphi\vert$, $\lambda = -1.254 + 2.400 i$};
	    \node[rotate=90] at (0.02, 0.55) {\footnotesize $\dot{x}$};
	    \node[rotate=0] at (0.44, 0.02) {\footnotesize $x$};
	\end{tikzonimage}\\
	\begin{tikzonimage}[trim=0 0 0 40, clip=true, width=0.45\textwidth]{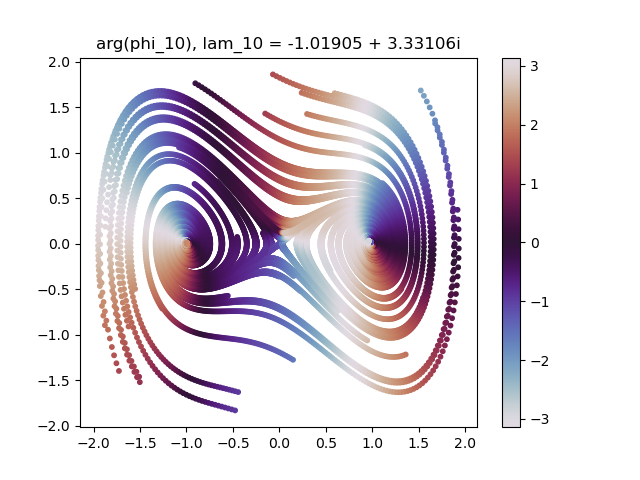}
	    \node[rotate=0] at (0.44, 1.075) {\footnotesize $\angle \varphi$, $\lambda = -1.019 + 3.331 i$};
	    \node[rotate=90] at (0.02, 0.55) {\footnotesize $\dot{x}$};
	    \node[rotate=0] at (0.44, 0.02) {\footnotesize $x$};
	\end{tikzonimage}
	\begin{tikzonimage}[trim=0 0 0 40, clip=true, width=0.45\textwidth]{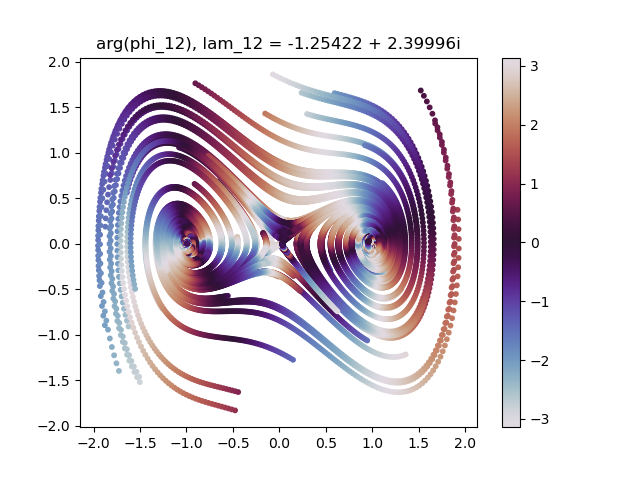}
	    \node[rotate=0] at (0.44, 1.075) {\footnotesize $\angle \varphi$, $\lambda = -1.254 + 2.400 i$};
	    \node[rotate=90] at (0.02, 0.55) {\footnotesize $\dot{x}$};
	    \node[rotate=0] at (0.44, 0.02) {\footnotesize $x$};
	\end{tikzonimage}
	\caption{Analog of \cref{fig:Duffing_eigenfunction0} for the eigenvalue closest to twice the eigenvalue of the linearization about the stable fixed points: EDMD (left); EM algorithm (right). }
	\label{fig:Duffing_eigenfunction2}
\end{figure}

We plot the corresponding approximate Koopman eigenfunctions obtained using both models in \cref{fig:Duffing_eigenfunction0}--\ref{fig:Duffing_eigenfunction2}.
The eigenfunction values for the HMM were computed using the state estimation-based approach described in \cref{sec:approximating_eigenfunction}.
We observe from \cref{fig:Duffing_eigenfunction0} that both methods identify an approximate eigenfunction that separates the two basins of attraction, but the EM-based HMM provides a more crisp boundary between the two basins.
This is possible because the EM-based approximation is not limited by the expressive power of the polynomial dictionary used in EDMD.
\Cref{fig:Duffing_eigenfunction1} shows the eigenfunction corresponding to the eigenvalue closest to that of the linearization about the stable fixed points. 
EDMD identifies these eigenvalues with less accuracy and the eigenfunction's angular behavior exhibits slight inaccuracies near the fixed point at $(1,0)$.
\Cref{fig:Duffing_eigenfunction2} shows the eigenfunction corresponding to twice the eigenvalues of the linearized system at the spiral fixed points.
We expect the angle to change by $4\pi$ around each stable fixed point.  This behavior is exhibited by the eigenfunctions approximated by the EM algorithm, but not by those approximated by EDMD.

The model we have learned for the Duffing equation using the EM algorithm still has a considerable amount of error both in predicting trajectories and in the higher eigenvalues and eigenfunctions.
These inaccuracies may be due to the fact that our model \cref{eqn:finite_dimensional_generator_IO_model} demands that the observations $\vect{y}$ are reconstructed as a linear combination of the latent variables $\vect{z}$. 
Hence, we are asking our learning procedure to trade off between finding an approximately Koopman-invariant latent space, and finding a latent space where the observed quantities can be linearly reconstructed.
These objectives are at odds when the observations demand a large superposition of Koopman eigenfunctions to be accurately represented, but the latent space is low-dimensional --- as is desirable for both learning and control purposes.
Numerical stability issues inherent with the explicit Euler discretization used to obtain \cref{eqn:main_HMM} also place limits on the spectra of the Koopman generator approximations and may prevent the model from taking advantage of large latent space dimensions.
Eigenvalues of the Koopman generator form a lattice in the complex plane, which can exit the stability region for \cref{eqn:main_HMM}.
However, this was not an issue here due to the small time step and correspondingly large region of stability for the explicit Euler scheme shown on the right of \cref{fig:Duffing_error_and_eigenvalues}.

As argued in \cite{Otto2019linearly}, it may be possible to significantly improve the performance without increasing the dimension of the latent state by also optimizing a nonlinear reconstruction map $\vect{y} = \vect{\tilde {g}}(\vect{z})$ that could be defined by a neural network.
In such a setup, one would have to modify the procedure in the E-step to account for the nonlinearity, $\vect{\tilde {g}}$, perhaps by using a sampling approach like unscented or particle filtering.

\subsection{Fluidic pinball with varying Reynolds number}
\label{subsec:FP}
\begin{figure}[bht!]
	\centering
	\begin{tikzonimage}[trim=0 0 0 0, clip=true, width=0.99\textwidth]{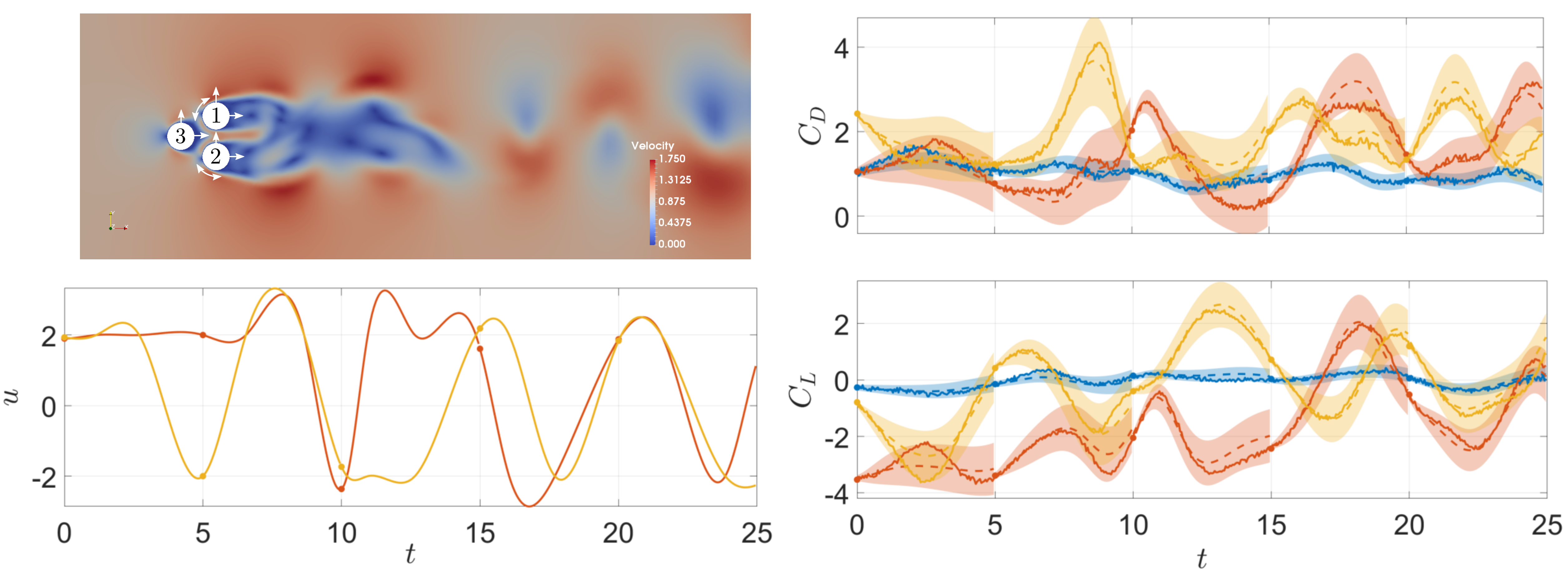}
	    \fill [white] (0.025,0.3) rectangle (0.0,0.36);
	    \node[rotate=90] at (0.015, 0.325) {\footnotesize input, $u$};
	    \fill [white] (0.25,0.0) rectangle (0.28,0.08);
	    \node[rotate=0] at (0.265, 0.05) {\footnotesize $t$};
	    \fill [white] (0.5,0.75) rectangle (0.53,0.83);
	    \node[rotate=90] at (0.515, 0.78) {\footnotesize drag, $C_D$};
	    \fill [white] (0.5,0.28) rectangle (0.53,0.38);
	    \node[rotate=90] at (0.515, 0.325) {\footnotesize lift, $C_L$};
	    \fill [white] (0.75,0.0) rectangle (0.78,0.08);
	    \node[rotate=0] at (0.765, 0.05) {\footnotesize $t$};
	    \node[rotate=0] at (0.12, 0.84) {\scriptsize \textcolor{white}{$u_1$}};
	    \node[rotate=0] at (0.12, 0.69) {\scriptsize \textcolor{white}{$u_2$}};
	\end{tikzonimage}
	\caption{Top left: Fluidic pinball setup. Other plots: Prediction (dashed lines) of lift $C_L$ and drag $C_D$ of the actuated fluidic pinball (cylinder 1: orange, cylinder 2: yellow, cylinder 3: blue) over five seconds, starting a different initial times denoted by the dots, and using the best model with $\dim{(\vect{z})}=25$ at $Re=225$, for which there was no data contained in the training data set. The shaded area denotes the $2\sigma$ confidence interval and the solid lines are the true trajectories.}
	\label{fig:FP_setup}
\end{figure}
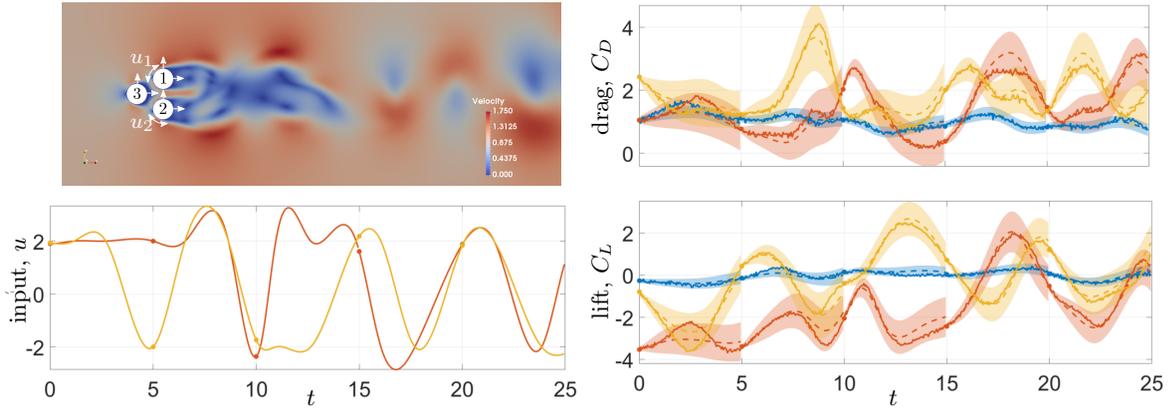

The final example is an unsteady flow around three cylinders arranged on an
equidistant triangle as shown in \cref{fig:FP_setup}, in a configuration known
as the \emph{fluidic pinball}.
The flow physics are governed by the two-dimensional incompressible
Navier--Stokes equations, and depend on a dimensionless parameter, the Reynolds number,
defined by $Re=U_{\infty}D/\nu$, where $U_{\infty}$ is the incoming velocity of the fluid, $\nu$
is the kinematic viscosity, and $D$ is
the diameter of each cylinder.
In our study, we nondimensionalize lengths by $D$, velocities by $U_{\infty}$, and time
by $D/U_{\infty}$.
In~\cite{DNMP19}, it was conjectured that the system exhibits chaotic behavior above the critical Reynolds number $Re\approx 115$. 
The chaotic regime has been studied in the control context before \cite{Peitz2020data,BPB+20}, and we are going to do so here as well by using our data-driven surrogate models for MPC.
Following \cite{BPB+20}, we provide input to the system by rotating cylinders $1$ and $2$ (see \cref{fig:FP_setup}) with counterclockwise angular velocities $u_1$ and $u_2$ measured in radians per dimensionless time.
For control, these inputs are constrained to the interval $u_1, u_2 \in [-2,2]$.

As the factor $1/Re$ enters linearly into the Navier--Stokes equations, we can simply treat it as an additional control input in our framework, i.e.,
\[
	\bar{u} = (u_1, u_2, 1/Re)
\]
and thus study the ability to cope with previously unseen parameter values.
To this end, we collect data at $Re=150$, $Re=200$ and $Re=250$, and study the model performance at $Re=175$, $Re=225$ and $Re=275$.

\begin{figure}
	\centering
	\begin{tikzonimage}[trim=0 0 0 0, clip=true, width=0.98\textwidth]{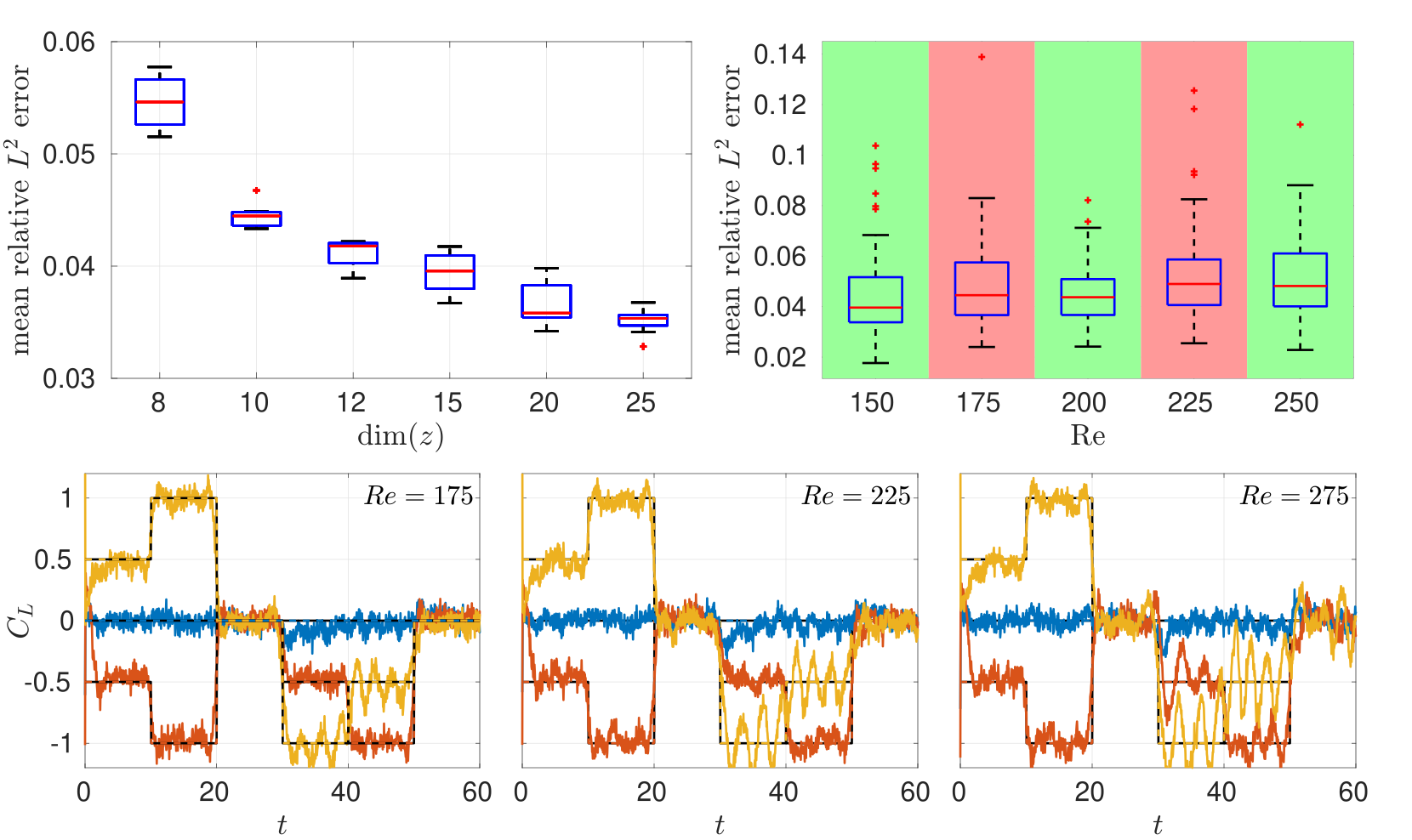}
	    \fill [white] (0.025,0.55) rectangle (0.0,0.925);
	    \node[rotate=90] at (0.015, 0.75) {\footnotesize mean relative $\ell^2$ error, $\epsilon$};
	    \fill [white] (0.25,0.46) rectangle (0.325,0.5);
	    \node[rotate=0] at (0.28, 0.475) {\footnotesize model dimension, $\dim(\vect{z})$};
	    \fill [white] (0.525,0.55) rectangle (0.5,0.925);
	    \node[rotate=90] at (0.515, 0.75) {\footnotesize mean relative $\ell^2$ error, $\epsilon$};
	    \fill [white] (0.75,0.46) rectangle (0.8,0.5);
	    \node[rotate=0] at (0.775, 0.475) {\footnotesize Reynolds number, $Re$};
	    \fill [white] (0.025,0.225) rectangle (0.0,0.3);
	    \node[rotate=90] at (0.005, 0.26) {\footnotesize lift coefficient, $C_L$};
	\end{tikzonimage}
	\caption{Top Left: Distribution of the training data prediction accuracy
          over 10 runs with different initial guesses $\mat{V}_i$ for various
          latent space dimensions. Top right: Box plot for the prediction
          accuracy of the best model ($\dim{(\vect{z})}=25$) over 500 test
          trajectories (100 for every $Re$), where only training data for
          $Re=150$, $200$ and $250$ was used. Here, ten observations are used to
          estimate the initial conditions for the latent state via the Kalman
          filter. 
          Bottom: Lift tracking
          for $Re=175$, $Re=225$ and $Re=275$ which were not contained in the
          training data.}
	\label{fig:FP_multiRe_boxplots_control}
\end{figure}

In the control setting, the task is to influence the (non-dimensional) drag $C_D$ or the (non-dimensional) lift forces $C_L$ on each cylinder by rotating the first and second cylinder without any direct observations of the flow field. 
The directions of the lift and drag forces on each cylinder are indicated by the horizontal and vertical arrows shown on each cylinder in the top left of \Cref{fig:FP_setup}.
We observe only the lift and drag on the three cylinders 
\[
    \vect{y} = (C_{D,1}, C_{D,2}, C_{D,3}, C_{L,1}, C_{L,2}, C_{L,3}),
\]
and collect data with $\Delta t=0.05$ from one long trajectory of length $500$ seconds per Reynolds number with a random smooth input $u(t)$ and measurement noise with variance $0.2$.
The input was generated by smooth spline interpolation on a random time grid with $0.5 < \Delta t < 2.5$ and with random values drawn from the square $[-3,3]^2$.
For the training, we then split these into 100 trajectories with 100 data points each.

In order to select the model dimension, we randomly initialize $10$ models at
each dimension $\dim(\vect{z}) = 8,10,12,15,20,$ and $25$, and train each until convergence when the likelihood fails to increase.
For each model, we use the optimal estimates of the initial states to forecast the training trajectories and compute the relative $\ell^2$ error
\begin{equation}\label{eq:L2Error}
	\epsilon = \frac{1}{500}\mathlarger{\sum}_{m=1}^{500} \sum_{l=1}^{100} \frac{\|\vect{\hat{y}}^{(m)}_l - \vect{y}^{(m)}_l\|_2}{\|\vect{y}^{(m)}_l\|_2}\Delta t.
\end{equation}
The distribution of these prediction errors among the randomly initialized models at each dimension are plotted in \cref{fig:FP_multiRe_boxplots_control}.
As we expect, models with larger dimensions make more accurate forecasts, and so we select the model with $\dim(\vect{z}) = 25$.
Over-fitting was not an issue here because the amount of data greatly exceeded
the number of free parameters in the largest model.

Representative predictions by the selected model on testing data are visualized in the plots on the right of \Cref{fig:FP_setup}.
Here, we use our learned model to make forecasts over each interval, where we use the last $10$ observations from the previous interval to provide a state estimate of the model's initial condition and its uncertainty.
We used a similar procedure to compute the prediction error on testing data consisting of $100$ trajectories at each Reynolds number in the top right plot of \Cref{fig:FP_multiRe_boxplots_control}.
We use the first $10$ samples along each trajectory to obtain a state estimate that is used to forecast the output at the remaining $90$ sample times and to compute the relative $\ell^2$ error.
We observe that the prediction error for the unseen Reynolds numbers $Re=175$
and $Re=225$, is not worse than for the Reynolds numbers contained in the
training set.

Finally, excellent MPC performance is observed for lift tracking at different
Reynolds numbers, as shown in the bottom plots of
\Cref{fig:FP_multiRe_boxplots_control}, where $Re=275$ even extends past the
Reynolds numbers used during training.
The only notable change is a performance decrease for the non-symmetric reference trajectories ($t=30$ to $t=50$), which increase the complexity of the dynamics. 


\section{Conclusion}
We have formulated the learning problem for matrix approximations of Koopman generators associated with nonlinear input-affine dynamical systems as a parameter estimation problem for a bilinear hidden Markov model.
This formulation allows us to identify the Koopman generator matrix approximations as well as covariance matrices quantifying the noise and model uncertainties from noisy actuated trajectories of the system's output using the Expectation-Maximization (EM) algorithm.
Using the EM algorithm overcomes the need in existing algorithms such as
Extended Dynamic Mode Decomposition (EDMD) to explicitly specify a dictionary of
functions in order to approximate the Koopman operator or generator.
Moreover, one doesn't need to access the system's full state; one can rely solely on trajectories of a few noisy, possibly nonlinear functions of the state that might be obtained from sensors in an experiment.

We demonstrated the EM-based learning technique and the performance of the
resulting models on several examples, in which we use the models to make
forecasts, estimate Koopman eigenfunctions, and implement model predictive control (MPC).
Our first example demonstrates that the technique is capable of learning the dynamics of a system whose evolution is described in a finite-dimensional Koopman-invariant subspace.
Using the Duffing equation in our second example, we show how the method can be used to estimate Koopman eigenfunctions without relying on a dictionary or other parametric families of functions like neural networks.
Finally, our fluidic pinball example demonstrates that the technique can be applied to forecast and control large-scale systems with partial observations.
Moreover, this example illustrates how we may treat parameters such as the Reynolds number as additional inputs in our framework.

Our analysis of the Duffing equation suggests several avenues for future work.
In this case, there is a finite-dimensional Koopman-invariant subspace, but an infinite number of Koopman eigenfunctions are needed to exactly reconstruct the system's output as a linear combination.
This leads to inaccuracies because the model formulation demands that the latent space is both Koopman-invariant and rich enough to linearly reconstruct the observations.
Our analysis after \cref{thm:maximization_step} and preliminary experiments suggest that the data requirements for the method grow rapidly with the model dimension.
We have suggested some ways to mitigate these issues via regularization, but more work is needed to elucidate the connections between the model dimension, data requirements, accuracy, and computational complexity of our proposed method.
Performance can likely also be improved by learning the parameters of a hidden Markov model where the observations are nonlinearly reconstructed via a neural network.
Another approach might leverage a model constructed via more accurate linear multistep methods employing implicit time discretization.

Our code implementing the EM algorithm for approximating Koopman generators is written in Python and is available at \url{https://github.com/samotto1/KoopmanGeneratorEM}.

\section*{Acknowledgements}
Calculations were performed by S.P.\ on resources provided by the Paderborn Center for Parallel Computing ($\mbox{PC}^2$).
This research was supported by the Army Research Office under grant number W911NF-17-1-0512 and the Air Force Office of Scientific Research under grant number FA9550-19-1-0005. 
S.E.O.\ was supported by the National Science Foundation Graduate Research Fellowship Program under Grant No. DGE-2039656. S.P.\ acknowledges support by the DFG Priority Programme 1962.


\bibliographystyle{siamplain}
\bibliography{jfull,References}


\appendix

\section{Maximization step: Proof of Theorem~\ref{thm:maximization_step}}
\label{app:maximization_step}
During maximization of the regularized evidence lower bound 
\begin{equation}
    \hat{L}_Q(\calP) - R(\calP)
    = \sum_{m=1}^M \mathbb{E}_{\mat{\hat{Z}}^{(m)}}\left[ \log{P_{\mat{Z}, \mat{Y}}(\mat{\hat{Z}}^{(m)}, \mat{Y}^{(m)};\ \calP)} \right] - \sum_{m=1}^M \mathbb{E}_{\mat{\hat{Z}}^{(m)}}\left[ \log{Q^{(m)}(\mat{\hat{Z}}^{(m)})} \right] - R(\calP)
\end{equation}
over the parameters $\mathcal{P}$ with fixed inference distributions $\{ Q^{(m)} \}_{m=1}^M$, we need only consider the term 
\begin{equation*}
    \sum_{m=1}^M \mathbb{E}_{\mat{\hat{Z}}^{(m)}}\left[ \log{P_{\mat{Z}, \mat{Y}}(\mat{\hat{Z}}^{(m)}, \mat{Y}^{(m)};\ \calP)} \right] - R(\calP),
\end{equation*}
since the remaining term $\sum_{m=1}^M \mathbb{E}_{\mat{\hat{Z}}^{(m)}}\left[ \log{Q^{(m)}(\mat{\hat{Z}}^{(m)})} \right]$ does not depend on $\calP$.
These $\calP$-dependent terms in the regularized ELBO decouple into three separate terms thanks to the Markov property of \cref{eqn:main_HMM}, as we show below in Proposition~\ref{prop:M_step_decoupling}.

Throughout, we denote the means and joint covariances of the inference distributions $Q^{(m)}$, $m=1,\ldots,M$ by
\begin{equation}
    \vmuhat_k^{(m)} = \mathbb{E}_{\mat{\hat{Z}}^{(m)}}\left[ \vect{\hat{z}}_k^{(m)} \right] 
    \qquad \text{and} \qquad
    \mSighat_{k,l}^{(m)} = \mathbb{E}_{\mat{\hat{Z}}^{(m)}}\left[ \big(\vect{\hat{z}}_k^{(m)} - \vmuhat_k^{(m)}\big) \big(\vect{\hat{z}}_l^{(m)} - \vmuhat_{l}^{(m)}\big)^T \right].
\end{equation}
The first decoupled $\calP$-dependent term in the regularized ELBO,
\begin{multline}
    L_1(\vmu_0, \mSig_0) := M\log\det{\left(2\pi \mSig_{0}\right)} \\
    + \Tr\Bigg\lbrace\mSig_{0}^{-1} \underbrace{\Bigg(\sum_{m=1}^M \Big[ \mSighat_{0,0}^{(m)} + \big(\vmuhat_0^{(m)} - \vmu_0\big)\big(\vmuhat_0^{(m)} - \vmu_0\big)^T \Big] + \gamma_0 \mat{I} \Bigg)}_{\mat{W}_1(\vmu_0)} \Bigg\rbrace,
    \label{eqn:IC_loss_fcn}
\end{multline}
is a loss function for the initial condition parameters.
The second term,
\begin{multline}
    L_2(\vect{c}_0, \mSig_{\vv}) := M(L+1)\log\det{\left(2\pi \mSig_{\vv}\right)} \\
    + \Tr\Bigg\lbrace \mSig_{\vv}^{-1} \underbrace{\Bigg( \sum_{m=1}^M \sum_{l=0}^L \Big[ \mat{\tilde{C}}\mSighat_{l,l}^{(m)}\mat{\tilde{C}}^T
    + \big(\vy_l^{(m)} - \vect{c}_0 - \mat{\tilde{C}}\vmuhat_l^{(m)}\big)\big(\vy_l^{(m)} - \vect{c}_0 - \mat{\tilde{C}}\vmuhat_l^{(m)}\big)^T \Big]
    + \gamma_{\vv} \mat{I} \Bigg)}_{\mat{W}_2(\vect{c}_0)} \Bigg\rbrace,
    \label{eqn:observation_loss_fcn}
\end{multline}
is a loss function for the observation map parameters.
Letting
\begin{multline}
    \mat{W}_3(\mat{\tilde{V}}_0,\ldots,\mat{\tilde{V}}_{\dim\vect{u}}) = \\
    \sum_{m=1}^M \sum_{l=0}^{L-1}\Big[ 
    \mSighat_{l+1,l+1}^{(m)} 
    - \mat{A}_l^{(m)}\mSighat_{l,l+1}^{(m)}
    - \mSighat_{l+1,l}^{(m)}(\mat{A}_l^{(m)})^T
    + \mat{A}_l^{(m)}\mSighat_{l,l}(\mat{A}_l^{(m)})^T \\
    + \big(\vmuhat_{l+1}^{(m)} - \mat{A}_l^{(m)}\vmuhat_l^{(m)} - \vect{b}_l^{(m)}\big)\big(\vmuhat_{l+1}^{(m)} - \mat{A}_l^{(m)}\vmuhat_l^{(m)} - \vect{b}_l^{(m)}\big)^T \Big]
    + \gamma_{\mat{G}} \Delta t^2 \sum_{i=0}^{\dim\vu}\mat{\tilde{V}}_i^T \mat{\tilde{V}}_i + \gamma_{\vw} \mat{I},
\end{multline}
the third term,
\begin{equation}
    L_3(\mat{\tilde{V}}_0,\ldots,\mat{\tilde{V}}_{\dim\vect{u}}, \mSig_{\vw}) = M L \log\det{\left(2\pi \mSig_{\vw}\right)}
    + \Tr\Big[ \mSig_{\vw}^{-1} \mat{W}_3(\mat{\tilde{V}}_0,\ldots,\mat{\tilde{V}}_{\dim\vect{u}}) \Big],
    \label{eqn:dynamics_loss_fcn}
\end{equation}
is a loss function for the dynamical parameters including the matrix approximations of the Koopman generators.
The precise decoupling of the regularized ELBO into the three loss functions is proved below.
\begin{proposition}[Decoupled objectives for maximization step]
\label{prop:M_step_decoupling}
The terms in the regularized ELBO depending on $\calP$ satisfy
\begin{multline}
    \sum_{m=1}^M \mathbb{E}_{\mat{\hat{Z}}^{(m)}}\left[ \log{P_{\mat{Z}, \mat{Y}}(\mat{\hat{Z}}^{(m)}, \mat{Y}^{(m)};\ \calP)} \right] - R(\calP) \\
    = -\frac{1}{2} \left[L_1(\vmu_0, \mSig_0) + L_2(\vect{c}_0, \mSig_{\vv}) + L_3(\mat{\tilde{V}}_0,\ldots, \mat{\tilde{V}}_{\dim\vect{u}}, \mSig_{\vw}) \right].
    \label{eqn:Expectation}
\end{multline}
\end{proposition}
\begin{proof}
We begin by considering a single trajectory with fixed $m$, removing it from the superscript, and take the sum over $m$ at the end.
To avoid cluttered equations for the time being, we also drop the explicit dependence of the probability distributions on the model parameters $\calP$ and remove the subscripts from probability densities where the appropriate random variables are obvious.
If we let $\vect{z}_{0:l} = (\vect{z}_0, \ldots, \vect{z}_{l})$ and $\vect{y}_{0:l} = (\vect{y}_{0}, \ldots, \vect{y}_l)$, then by the Markov property we have
\begin{equation}
\begin{split}
    P(\mat{Z}, \mat{Y}) &= P(\vect{z}_{0:L}, \vect{y}_{0:L}) \\
    &= P(\vect{z}_{0:L-1}, \vect{y}_{0:L-1}) \cdot P\left(\vz_L, \vy_L\ \vert\ \vect{z}_{0:L-1}, \vect{y}_{0:L-1}\right) \\
    &= P(\vect{z}_{0:L-1}, \vect{y}_{0:L-1}) \cdot P\left(\vz_L, \vy_L\ \vert\ \vz_{L-1}\right) \\
    &\vdots \\
    &= P\left(\vz_0, \vy_0\right)\cdot \prod_{l=0}^{L-1} P\left(\vz_{l+1}, \vy_{l+1}\ \vert\ \vz_{l}\right).
\end{split}
\end{equation}
Making use of the conditional independence of the observation $\vy_{l+1}$ and previous state $\vz_{l}$ given $\vz_{l+1}$, we obtain
\begin{equation}
    P(\mat{Z}, \mat{Y}) = P\left( \vz_0 \right) \cdot
    \prod_{l=0}^{L} P\left(\vy_{l} \ \vert \ \vz_{l}\right) \cdot 
    \prod_{l=0}^{L-1} P\left(\vz_{l+1}\ \vert\ \vz_{l}\right).
\end{equation}
Recalling our dynamical model \cref{eqn:main_HMM}, the log joint probability is given by
\begin{equation}
\begin{split}
    \log{P(\mat{Z}, \mat{Y})} =\ & - \frac{1}{2}\log\det{\left(2\pi \mSig_{0}\right)} - \frac{1}{2}\left(\vz_0 - \vmu_0\right)^T\mSig_0^{-1}\left(\vz_0 - \vmu_0\right) \\
    & -\frac{1}{2}(L+1)\log\det{\left(2\pi \mSig_{\vv}\right)} - \frac{1}{2}\sum_{l=0}^L\left(\vy_l - \vect{c}_0 - \mat{\tilde{C}}\vz_l \right)^T\mSig_{\vv}^{-1}\left(\vy_l - \vect{c}_0 - \mat{\tilde{C}}\vz_l \right) \\
    & -\frac{1}{2}L\log\det{\left(2\pi \mSig_{\vw}\right)}
    - \frac{1}{2}\sum_{l=0}^{L-1}\left(\vz_{l+1} - \mat{A}_l\vz_l - \vect{b}_l \right)^T\mSig_{\vw}^{-1}\left(\vz_{l+1} - \mat{A}_l\vz_l - \vect{b}_l \right).
\end{split}
\end{equation}

Taking the expectation with respect to the inference distribution, we obtain
\begin{equation}
    \mathbb{E}_{\mat{\hat{Z}}}\left[\log{P_{\mat{Z},\mat{Y}}(\mat{\hat{Z}}, \mat{Y};\ \calP)}\right] = -\frac{1}{2} \left[ \tilde{L}_1(\vmu_0, \mSig_0) + \tilde{L}_2(\mat{\tilde{C}}, \vect{c}_0, \mSig_{\vv}) + \tilde{L}_3(\mat{\tilde{V}}_0,\ldots,\mat{\tilde{V}}_q, \mSig_{\vw})\right],
\end{equation}
where
\begin{equation}
    \tilde{L}_1(\vmu_0, \mSig_0) = \log\det{\left(2\pi \mSig_{0}\right)} + \Tr\left[\mSig_{0}^{-1}\left( \mSighat_{0,0} + \left(\vmuhat_0 - \vmu_0\right)\left(\vmuhat_0 - \vmu_0\right)^T \right) \right]
\end{equation}
\begin{multline}
    \tilde{L}_2(\mat{\tilde{C}}, \vect{c}_0, \mSig_{\vv}) = (L+1)\log\det{\left(2\pi \mSig_{\vv}\right)} \\
    + \Tr\Bigg\lbrace \mSig_{\vv}^{-1} \sum_{k=0}^L \Big[ \mat{\tilde{C}}\mSighat_{k,k}\mat{\tilde{C}}^T
    + \left(\vy_k - \vect{c}_0 - \mat{\tilde{C}}\vmuhat_k\right)\left(\vy_k - \vect{c}_0 - \mat{\tilde{C}}\vmuhat_k\right)^T \Big] \Bigg\rbrace
\end{multline}
\begin{multline}
    \tilde{L}_3(\mat{\tilde{V}}_0,\ldots,\mat{\tilde{V}}_q, \mSig_{\vw}) = L \log\det{\left(2\pi \mSig_{\vw}\right)}
    + \Tr\Bigg\lbrace \mSig_{\vw}^{-1} \sum_{k=0}^{L-1}\Big[ 
    \mSighat_{k+1,k+1} 
    - \mat{A}_k\mSighat_{k,k+1} \\
    - \mSighat_{k+1,k}\mat{A}_k^T
    + \mat{A}_k\mSighat_{k,k}\mat{A}_k^T
    + \left(\vmuhat_{k+1} - \mat{A}_k\vmuhat_k - \vect{b}_k\right)\left(\vmuhat_{k+1} - \mat{A}_k\vmuhat_k - \vect{b}_k\right)^T \Big] \Bigg\rbrace
\end{multline}
The final result is obtained by summing over $m$ and subtracting the regularization term $R(\calP)$ defined by \cref{eqn:regularization_function}.
\end{proof}

We observe that by Proposition~\ref{prop:M_step_decoupling}, it suffices to minimize the three terms $L_1$, $L_2$, and $L_3$ separately, which each have the same form,
\begin{equation}
    L_i(\mSig_i, \mathcal{P}_i) = \alpha_i \log{\det{(2\pi \mSig_i)}} + \Tr{\left[ \mSig_i^{-1} \mat{W}_i(\mathcal{P}_i) \right]},
\end{equation}
where $\mSig_i$ is a covariance matrix to be determined, $\alpha_i > 0$ is a constant, and $\mat{W}_i$ is a symmetric, positive semidefinite matrix-valued function of the remaining parameters $\mathcal{P}_i\subset\mathcal{P}\setminus \{\mSig_i\}$ to be optimized.
In particular, after some algebraic manipulation, these functions are given by the following quadratic forms:
\begin{equation}
    \mat{W}_1(\vmu_0) 
    = \underbrace{\gamma_0 \mat{I} + \sum_{m=1}^M \mSighat_{0,0}^{(m)}}_{\mat{S}_1} 
    - \vmu_0 \left( \sum_{m=1}^M \vmuhat_0^{(m)} \right)^T
    - \left( \sum_{m=1}^M \vmuhat_0^{(m)} \right) \vmu_0^T
    + M \vmu_0 \vmu_0^T,
    \label{eqn:W1}
\end{equation}
\begin{multline}
    \mat{W}_2(\vect{c}_0)
    = \underbrace{\gamma_{\vv} \mat{I}
    + \sum_{m=1}^M \sum_{l=0}^L \mat{\tilde{C}}\mSighat_{l,l}^{(m)}\mat{\tilde{C}}^T}_{\mat{S}_2}
    - \vect{c}_0 \left[ \sum_{m=1}^M \sum_{l=0}^L \big( \vy_l^{(m)} - \mat{\tilde{C}}\vmuhat_l^{(m)} \big) \right]^T \\
    - \left[ \sum_{m=1}^M \sum_{l=0}^L \big( \vy_l^{(m)} - \mat{\tilde{C}}\vmuhat_l^{(m)} \big) \right] \vect{c}_0^T
    + M (L+1) \vect{c}_0 \vect{c}_0^T,
    \label{eqn:W2}
\end{multline}
and
\begin{equation}
    \mat{W}_3(\mat{\tilde{V}})
    =  \mat{S}_3 + \Delta t^2 \left( - \mat{H} \mat{\tilde{V}} - \mat{\tilde{V}}^T \mat{H}^T + \mat{\tilde{V}}^T \mat{G} \mat{\tilde{V}}\right),
    \qquad \mbox{where} \qquad \mat{\tilde{V}} = \begin{bmatrix} \mat{\tilde{V}}_0 \\ \vdots \\ \mat{\tilde{V}}_{\dim\vect{u}} \end{bmatrix},
    \label{eqn:W3}
\end{equation}
and
\begin{equation}
    \mat{S}_3 = \gamma_{\vw}\mat{I} + \sum_{m=1}^M \sum_{l=0}^{L-1} \Big[ \mSighat_{l+1,l+1}^{(m)} - \mSighat_{l,l+1}^{(m)} - \mSighat_{l+1,l}^{(m)} + \mSighat_{l,l}^{(m)}
    + \big( \vmuhat_{l+1}^{(m)} - \vmuhat_{l}^{(m)} \big) \big( \vmuhat_{l+1}^{(m)} - \vmuhat_{l}^{(m)} \big)^T \Big].
\end{equation}

We use \cref{lem:form_of_maximization_objective}, below, to show that each function $L_i(\mSig_i, \mathcal{P}_i)$ has a unique minimizer that can be computed explicitly.
\begin{lemma}
    \label{lem:form_of_maximization_objective}
    Let $\alpha > 0$ and suppose that $\calP'\mapsto\mat{W}(\calP')$ is a continuously differentiable positive semidefinite $n\times n$ matrix-valued function of finitely many real parameters $\calP'$ satisfying $\Tr(\mat{W}(\calP')) \to \infty$ as $\Vert \calP'\Vert \to \infty$ with respect to a norm $\Vert \cdot \Vert$.
    Suppose that there is a set of parameters $\calP'_*$ such that $\mat{W}(\calP'_*)$ is positive-definite and $\calP'_*$ is the unique solution of 
    $\frac{\partial}{\partial \calP'}\Tr(\mSig^{-1} \mat{W}(\calP')) = 0$ for every positive-definite matrix $\mSig$.
    Then $\big(\frac{1}{\alpha}\mat{W}(\calP'_*),\ \calP'_*\big)$ is the unique minimizer of
    \begin{equation}
        L(\mSig, \mathcal{P}') = \alpha \log\det\left(2\pi \mSig\right) + \Tr\left[\mSig^{-1} \mat{W}(\mathcal{P}') \right]
    \end{equation}
    over the set where $\mSig$ is positive-definite.
\end{lemma}
\begin{proof}
    With $\mSig$ being a fixed positive-definite matrix, we begin by showing that the function $L_{\mSig}: \calP' \mapsto L(\mSig, \calP')$ attains its minimum.
    Suppose that $\{\calP'_k\}_{k=1}^{\infty}$ is a sequence such that
    \begin{equation}
        \lim_{k\to\infty} L_{\mSig}(\calP'_k)
        = \inf_{\calP'} L_{\mSig}(\calP'),
    \end{equation}
    where the infemum is either finite or $-\infty$.
    First, it is clear that $\{\calP'_k\}_{k=1}^{\infty}$ is bounded, for if not, there is a subsequence with $\Vert \calP'_{k_m}\Vert \to \infty$ as $m\to\infty$.
    This implies that
    \begin{equation}
        L_{\mSig}(\calP'_{k_m}) 
        \geq \alpha \log\det\left(2\pi \mSig\right) + \frac{1}{\Vert \mSig \Vert} \Tr(\mat{W}(\calP'_{k_m})) \to \infty, 
    \end{equation}
    which is a contradiction.
    Since $\{\calP'_k\}_{k=1}^{\infty}$ is bounded in a finite-dimensional real vector space, there is a convergent subsequence $\calP'_{k_m} \to \calP'_{0}$.
    By continuity of $L_{\mSig}$, the limiting element minimizes $L_{\mSig}$ since
    \begin{equation}
        L_{\mSig}(\calP'_0) = \lim_{k\to\infty} L_{\mSig}(\calP'_k)
        = \inf_{\calP'} L_{\mSig}(\calP').
    \end{equation}
    Differentiating $L_{\mSig}$, we find that $\calP'_{0}$ satisfies $\frac{\partial \Tr(\mSig^{-1} \mat{W})}{\partial \calP'}(\calP'_0) = 0$.
    By our assumption it follows that $\calP'_0 = \calP'_*$.
    Since every convergent subsequence of the bounded sequence $\{ \calP'_k \}_{k=1}^{\infty}$ converges to $\calP'_*$ by the above argument, it follows that the entire sequence converges $\calP'_k \to \calP'_*$.

    Now suppose that $\{(\mSig_k, \calP'_k)\}$ is a sequence satisfying
    \begin{equation}
        \lim_{k\to\infty} L(\mSig_k, \calP'_k) 
        = \inf_{\mSig \succ 0, \calP'} L(\mSig, \calP').
    \end{equation}
    Since $L(\mSig_k, \calP'_k) \geq L(\mSig_k, \calP'_*)$, we also have
    \begin{equation}
        \lim_{k\to\infty} L(\mSig_k, \calP'_*) = \inf_{\mSig \succ 0, \calP'} L(\mSig, \calP').
    \end{equation}
    Let $c > 0$ denote the smallest eigenvalue of $\mat{W}(\calP'_*)$ and let $\lambda_1(\mSig_k) \geq \cdots \geq \lambda_n(\mSig_k) > 0$ denote the eigenvalues of $\mSig_k$ arranged in non-increasing order.
    With these definitions, we have
    \begin{equation}
        L(\mSig_k, \calP'_*) 
        \geq \sum_{i=1}^{n}\left[ \alpha \log \big(2\pi \lambda_i(\mSig_k)\big) + \frac{c}{\lambda_i(\mSig_k)} \right].
        \label{eqn:eigenvalue_bound_for_general_log_likelihood_component}
    \end{equation}
    Since the matrices $\mSig_k$ are symmetric and positive-definite, we have $\lambda_1(\mSig_k) = \Vert \mSig_k\Vert$.
    This implies that
    \begin{equation}
        L(\mSig_k, \calP'_*) \geq \alpha \log\big( 2\pi \Vert \mSig_k\Vert \big),
    \end{equation}
    meaning that $\{ \mSig_k \}_{k=1}^{\infty}$ is a bounded sequence.
    Passing to a convergent subsequence, we can assume that $\mSig_{k_m} \to \mSig_*$ as $m\to\infty$, where $\mSig_*$ is positive semi-definite.
    
    To prove that $\mSig_*$ is positive-definite, we first observe that $\lambda_i(\mSig_{k_m}) \to \lambda_i(\mSig_*)$ 
    thanks to Corollary~4.3.15 in \cite{Horn2013matrix}, which implies that
    \begin{equation}
        \vert \lambda_i(\mSig_{k_m}) - \lambda_i(\mSig_*) \vert \leq \Vert \mSig_{k_m} - \mSig_* \Vert \to 0.
    \end{equation}
    Using the concavity of $\log$, it can be shown that
    \begin{equation}
        \alpha \log(2\pi \lambda) + \frac{c}{\lambda}
        \geq \alpha \big(\log(2\pi\beta) + 1 \big) + \frac{c - \alpha \beta}{\lambda}
    \end{equation}
    for every $\lambda, \beta > 0$.
    Taking $\beta = c/\alpha$ yields the result that $\alpha \log(2\pi \lambda) + \frac{c}{\lambda} \geq \alpha (\log(2\pi c/\alpha) + 1)$ is bounded below.
    Taking $\beta = c/(2\alpha)$, yields
    \begin{equation}
        \alpha \log(2\pi \lambda) + \frac{c}{\lambda} 
        \geq \alpha \big(\log(\pi c / \alpha) + 1 \big) + \frac{c}{2 \lambda} \to \infty
    \end{equation}
    as $\lambda \to 0$.
    Using these results in \cref{eqn:eigenvalue_bound_for_general_log_likelihood_component} yields
    \begin{equation}
        L(\mSig_{k_m}, \calP'_*) 
        \geq (n-1)\alpha \big(\log(2\pi c/\alpha) + 1 \big) 
        + \alpha \big(\log(\pi c / \alpha) + 1 \big) + \frac{c}{2 \lambda_{n}(\mSig_{k_m})}
        \to \infty
    \end{equation}
    if $\lambda_{n}(\mSig_*) = 0$.
    Since $L(\mSig_{k_m}, \calP'_*) \to \infty$ is a contradiction, we must have $\lambda_{n}(\mSig_*) > 0$, meaning that $\mSig_*$ is positive-definite.
    By continuity of $L$, it follows that
    \begin{equation}
        L(\mSig_*, \calP'_*) 
        = \lim_{m\to\infty} L(\mSig_{k_m}, \calP'_*) 
        = \inf_{\mSig \succ 0, \calP'} L(\mSig, \calP').
    \end{equation}

    Differentiating $L$ with respect to $\mSig$ at $(\mSig_*, \calP'_*)$ yields
    \begin{equation}
    \begin{split}
        0 = \frac{\partial L}{\partial \mSig} (\mSig_*, \calP'_*) \delta\mSig 
        &= \alpha\Tr\left( \mSig_*^{-1}\delta\mSig \right) - \Tr\left[ \mSig_*^{-1}(\delta\mSig)\mSig_*^{-1} \mat{W}(\mathcal{P}'_*) \right] \\
        &= \Tr\left[ \left(\alpha\mSig_*^{-1} - \mSig_*^{-1} \mat{W}(\mathcal{P}'_*) \mSig_*^{-1} \right) \delta\mSig \right],
    \end{split}
    \end{equation}
    for every $n\times n$ matrix $\delta\mSig$.
    Here we have used the formulas for differentiating the log determinant and the matrix inverse given by Theorems~2~and~3 in Section~8.4 of \cite{Magnus2007matrix}.
    Since the above expression holds for arbitrary $\delta\mSig$, we must have
    \begin{equation}
        \alpha \mSig_*^{-1} = \mSig_*^{-1} \mat{W}(\mathcal{P}'_*) \mSig_*^{-1},
    \end{equation}
    for otherwise we could choose $\delta\mSig = \left[\alpha \mSig_*^{-1} - \mSig_*^{-1} \mat{W}(\mathcal{P}'_*) \mSig_*^{-1}\right]^T$ and produce a contradiction.
    Multiplying on both sides by $\mSig_*$ and dividing by $\alpha$ yields
    \begin{equation}
        \mSig_* = \frac{1}{\alpha} \mat{W}(\calP'_*).
    \end{equation}
    Since this must be the limit of any converging subsequence $\{\mSig_{k_m}\}_{m=1}^{\infty}$ of the bounded sequence $\{\mSig_k\}_{k=1}^{\infty}$, it follows that the entire sequence converges to $\mSig_k \to \mSig_* = \frac{1}{\alpha} \mat{W}(\calP'_*)$.

    To prove uniqueness, suppose that $(\mSig_{**}, \calP'_{**})$ is another minimizer of $L$ and let $(\mSig_k, \calP'_k) \to (\mSig_{**}, \calP'_{**})$.
    For sufficiently large $k$, $\mSig_k$ is positive-definite, and it follows by continuity of $L$ that $L(\mSig_k, \calP'_k) \to L(\mSig_*, \calP'_*)$.
    By the above argument, we have $\mSig_k \to \mSig_*$, meaning that $\mSig_{**} = \mSig_{*}$.
    Again, by continuity of $L$ we have 
    \begin{equation}
        \lim_{k\to\infty} L_{\mSig_*}(\calP'_k) 
        = L(\mSig_{**}, \calP'_{**})
        = \inf_{\mSig \succ 0, \calP'} L(\mSig, \calP')
        = \inf_{\calP'} L_{\mSig_*}(\calP').
    \end{equation}
    By uniqueness of the minimizer of $L_{\mSig_*}$, it follows that $\calP'_k \to \calP'_*$, meaning that $\calP'_{**} = \calP'_*$.
\end{proof}

It is easy to verify the hypotheses of the lemma for each $\mat{W}_i$ under the positive-definiteness assumptions of \cref{thm:maximization_step}.
Considering $\mat{W}_1$, we have
\begin{equation}
    \Tr\big[ \mSig^{-1} \mat{W}_1(\vmu_0) \big]
    = \Tr\big( \mSig^{-1} \mat{S}_1 \big) 
    - 2 \left( \sum_{m=1}^{M} \vmuhat_0^{(m)} \right) \mSig^{-1} \vmu_0
    + M \vmu_0^T \mSig^{-1} \vmu_0.
\end{equation}
Evidently, $\Tr[ \mat{W}_1(\vmu_0) ] \to \infty$ when $\Vert \vmu_0 \Vert \to \infty$.
Differentiating yields
\begin{equation}
    \frac{\partial}{\partial \vmu_0} \Tr\big[ \mSig^{-1} \mat{W}_1(\vmu_0) \big] \delta \vmu_0
    = 2 M \vmu_0^T \mSig^{-1} \delta\vmu_0 - 2 \left( \sum_{m=1}^{M} \vmuhat_0^{(m)} \right) \mSig^{-1} \delta\vmu_0
\end{equation}
for every $\delta \vmu_0$.
Therefore, $\frac{\partial \Tr(\mSig^{-1} \mat{W}_1)}{\partial \vmu_0} (\left.\vmu_0\right._*) = 0$ if and only if
\begin{equation}
    \left.\vmu_0\right._* = \frac{1}{M} \sum_{m=1}^{M} \vmuhat_0^{(m)}.
\end{equation}
Since \cref{eqn:initial_condition_covariance} is equal to $\frac{1}{M}\mat{W}_1(\left.\vmu_0\right._*)$, our assumption that \cref{eqn:initial_condition_covariance} is positive-definite is equivalent to the assumption that $\mat{W}_1(\left.\vmu_0\right._*)$ is positive-definite.
The argument for $\mat{W}_2$ is analogous.
Therefore, the conditions of \cref{lem:form_of_maximization_objective} hold for $L_1$ and $L_2$.

Finally, considering $\mat{W}_3$, we obtain
\begin{equation}
    \Tr\big[ \mSig^{-1} \mat{W}_3(\mat{\tilde{V}}) \big]
    = \Tr\big( \mSig^{-1} \mat{S}_3 \big)
    - 2 \Delta t^2 \Tr\big( \mSig^{-1} \mat{H} \mat{\tilde{V}} \big)
    + \Delta t^2 \Tr\big( \mSig^{-1} \mat{\tilde{V}}^T \mat{G} \mat{\tilde{V}} \big),
\end{equation}
where we have used the symmetry of $\mSig^{-1}$ and the permutation identity for the trace to write $\Tr(\mSig^{-1} \mat{\tilde{V}}^T \mat{H}^T ) = \Tr( \mat{H} \mat{\tilde{V}} \mSig^{-1} ) = \Tr( \mSig^{-1} \mat{H} \mat{\tilde{V}})$.
Thanks to our assumption that $\mat{G}$ is positive-definite, we evidently have $\Tr[ \mat{W}_3(\mat{\tilde{V}}) ] \to \infty$ when $\Vert \mat{\tilde{V}} \Vert \to \infty$.
Differentiating yields
\begin{equation}
    \frac{\partial}{\partial \mat{\tilde{V}}} \Tr\big[ \mSig^{-1} \mat{W}_3(\mat{\tilde{V}}) \big] \delta \mat{\tilde{V}}
    = 2\Delta t \Tr\left[ \mSig^{-1} \big( \mat{\tilde{V}}^T \mat{G} - \mat{H} \big) \delta\mat{\tilde{V}} \right]
    \qquad \forall \delta\mat{\tilde{V}},
\end{equation}
where have again used symmetry of $\mSig^{-1}$, symmetry of $\mat{G}$, and the permutation identity for the trace to express $\Tr(\mSig^{-1} \delta \mat{\tilde{V}}^T \mat{G} \mat{\tilde{V}} ) = \Tr( \mat{\tilde{V}}^T \mat{G} \delta \mat{\tilde{V}} \mSig^{-1} ) = \Tr( \mSig^{-1} \mat{\tilde{V}}^T \mat{G} \delta \mat{\tilde{V}} )$.
Therefore, we have $\frac{\partial \Tr(\mSig^{-1} \mat{W}_3 ) }{\partial \mat{\tilde{V}}}(\mat{\tilde{V}}_*) = 0$ if and only if
\begin{equation}
    \mat{\tilde{V}}_* = \mat{G}^{-1} \mat{H}^{T}.
\end{equation}
Since \cref{eqn:process_noise_covariance} is equal to $\frac{1}{M L} \mat{W}_3(\mat{\tilde{V}}_*)$, our assumption that \cref{eqn:process_noise_covariance} is positive-definite is equivalent to assuming that $\mat{W}_3(\mat{\tilde{V}}_*)$ is positive-definite.
Hence, the conditions of \cref{lem:form_of_maximization_objective} are verified for $L_3$.
Using the lemma to minimize $L_1$, $L_2$, and $L_3$ completes the proof.

\section{Expectation step: deriving the filtering and smoothing equations}
\label{app:expectation_step}

We present a complete derivation of the equations presented in \cref{subsec:Estep} by closely following the derivation by Yu et al.\ \cite{Yu2004derivation}.
Recall that the hidden state variables $\vect{z}_l$ along a trajectory are governed by \cref{eqn:E_step_dynamical_model}, that is
\begin{equation}
\begin{split}
    \vect{z}_{l+1} &= \mat{A}_l\vect{z}_l + \vect{b}_l + \vect{w}_l \\
    \vect{y}_l &= \vect{c}_0 + \mat{\tilde{C}}\vect{z}_l + \vect{v}_l,
\end{split}
\label{eqn:E_step_dynamical_model_app}
\end{equation}
where $\mat{A}_l$, $\vect{b}_l$, $\vect{c}_0$, and $\mat{\tilde{C}}$ are known. 
The initial condition, process noise, and measurement noise have independent Gaussian distributions
\begin{equation}
    \vect{z}_0 \sim \mathcal{N}(\vmu_0, \mSig_0), \qquad
    \vect{w}_l \sim \mathcal{N}(\vect{0}, \mSig_{\vect{w}}), \qquad
    \vect{v}_l \sim \mathcal{N}(\vect{0}, \mSig_{\vect{v}}),
\end{equation}
with known means and covariances.
Denoting $\mat{A}_{l,k} = \mat{A}_{l-1} \mat{A}_{l-2} \cdots \mat{A}_k$ for $l> k$ and $A_{k,k} = I$, the trajectory of \cref{eqn:E_step_dynamical_model_app} is given explicitly by
\begin{equation}
\begin{split}
    \vect{z}_{l} &= \mat{A}_{l,0} \vect{z}_0 + \sum_{k=0}^{l-1} \mat{A}_{l,k+1}(\vect{b}_{k} + \vw_{k}) \\
    \vect{y}_{l} &= \vect{c}_0 + \vect{v}_l + \mat{\tilde{C}}\mat{A}_{l,0} \vect{z}_0 + \sum_{k=0}^{l-1}\mat{\tilde{C}}\mat{A}_{l,k+1}(\vect{b}_{k} + \vw_{k}).
\end{split}
\end{equation}
We also denote $\vect{z}_{k:l} = (\vect{z}_k, \vect{z}_{k+1}, \ldots, \vect{z}_l)$, $k \leq l$ and likewise for other time-dependent quantities.
Since the vector $(\vect{z}_{0:L}, \vect{y}_{0:L})$ is a linear function of the Gaussian random variables $\vect{z}_0$, $\vect{w}_{0:L-1}$, and $\vect{v}_{0:L}$, it follows that $(\vect{z}_{0:L}, \vect{y}_{0:L})$ has a Gaussian distribution.
In particular, the conditional distributions $P(\vect{z}_l\vert\vect{y}_{0:L})$ and $P(\vect{z}_{l:l+1}\vert\vect{y}_{0:L})$ are Gaussian.
We seek to compute the conditional mean
\begin{equation}
    \vmuhat_l = \mathbb{E}\big[ \vect{z}_l \ \vert \ \vect{y}_{0:L} \big]
\end{equation}
and covariance matrices
\begin{equation}
    \mSighat_{l,l} = \mathbb{E}\big[ (\vect{z}_l - \vmuhat_l) (\vect{z}_l - \vmuhat_l)^T \ \vert \ \vect{y}_{0:L} \big], \qquad
    \mSighat_{l,l+1} = \mathbb{E}\big[ (\vect{z}_l - \vmuhat_l) (\vect{z}_{l+1} - \vmuhat_{l+1})^T \ \vert \ \vect{y}_{0:L} \big],
\end{equation}
which are parameters of the Gaussian distributions $P(\vect{z}_l\vert\vect{y}_{0:L})$ and $P(\vect{z}_{l:l+1}\vert\vect{y}_{0:L})$.

\subsection{Forward recursion (filtering)}
We begin by deriving the forward recursion to compute the conditional mean 
\begin{equation}
    \vmuhat_{l\vert l} = \mathbb{E}\big[ \vect{z}_l \ \vert \ \vect{y}_{0:l} \big]
\end{equation}
and conditional covariance matrices
\begin{equation}
    \mSighat_{l,l\vert l} = \mathbb{E}\big[ (\vect{z}_l - \vmuhat_l) (\vect{z}_l - \vmuhat_l)^T \ \vert \ \vect{y}_{0:l} \big], \qquad
    \mSighat_{l,l\vert l-1} = \mathbb{E}\big[ (\vect{z}_l - \vmuhat_l) (\vect{z}_l - \vmuhat_l)^T \ \vert \ \vect{y}_{0:l-1} \big].
\end{equation}
Using Bayes' theorem, and the Markov property of \cref{eqn:E_step_dynamical_model_app} we have
\begin{equation}
\begin{aligned}
    \log P(\vect{z}_l\vert \vect{y}_{0:l}) &=
    \log P(\vect{y}_l \vert \vect{z}_l, \vect{y}_{0:l-1}) 
    + \log P(\vect{z}_l \vert \vect{y}_{0:l-1})
    - \log P(\vect{y}_l\vert \vect{y}_{0:l-1}) \\
    &=
    \log P(\vect{y}_l \vert \vect{z}_l) 
    + \log P(\vect{z}_l \vert \vect{y}_{0:l-1})
    - \log P(\vect{y}_l\vert \vect{y}_{0:l-1})
\end{aligned}
\end{equation}
for every $l > 0$.
To determine $\vmuhat_{l\vert l}$ and $\mSighat_{l,l\vert l}$, it suffices to identify the terms of
\begin{equation}
\begin{aligned}
    \log P(\vect{z}_l\vert \vect{y}_{0:l}) 
    &= -\frac{1}{2} (\vect{z}_l - \vmuhat_{l\vert l})^T \mSighat_{l,l\vert l}^{-1} (\vect{z}_l - \vmuhat_{l\vert l}) + \cdots \\
    &= -\frac{1}{2} \vect{z}_l^T \mSighat_{l,l\vert l}^{-1} \vect{z}_l + \vect{z}_l^T \mSighat_{l,l\vert l}^{-1} \vmuhat_{l\vert l} + \cdots
\end{aligned}
\end{equation}
with linear and quadratic dependence on $\vect{z}_l$.
For short, terms in a sum that do not depend on $\vect{z}_l$ will be denoted with an ellipis ($\cdots$).
First, we note that $\log P(\vect{y}_l\vert \vect{y}_{0:l-1})$ does not depend on $\vect{z}_l$ at all.
By the definition of \cref{eqn:E_step_dynamical_model_app}, we have
\begin{equation}
\begin{aligned}
    \log P(\vect{y}_l \vert \vect{z}_l)
    &= -\frac{1}{2}( \vect{y}_l - \vect{c}_0 - \mat{\tilde{C}}\vect{z}_l )^T \mSig_{\vv}^{-1} ( \vect{y}_l - \vect{c}_0 - \mat{\tilde{C}}\vect{z}_l ) + \cdots \\
    &= -\frac{1}{2} \vect{z}_l^T \mat{\tilde{C}}^T \mSig_{\vv}^{-1} \mat{\tilde{C}}\vect{z}_l + \vect{z}_l^T \mat{\tilde{C}}^T \mSig_{\vv}^{-1} (\vect{y}_l - \vect{c}_0) + \cdots.
\end{aligned}
\end{equation}
We also have
\begin{equation}
    \log P(\vect{z}_l \vert \vect{y}_{0:l-1}) 
    = -\frac{1}{2} \vect{z}_l^T \mSighat_{l,l\vert l-1}^{-1} \vect{z}_l + \vect{z}_l^T \mSighat_{l,l\vert l-1}^{-1} \vmuhat_{l\vert l-1} + \cdots.
\end{equation}
Putting these expressions together yields
\begin{equation}
    \mSighat_{l,l\vert l}^{-1} = \mat{\tilde{C}}^T \mSig_{\vv}^{-1} \mat{\tilde{C}} + \mSighat_{l,l\vert l-1}^{-1}, 
    \quad \mbox{and} \quad
    \vmuhat_{l\vert l} = \mSighat_{l,l\vert l} \mat{\tilde{C}}^T \mSig_{\vv}^{-1} (\vect{y}_l - \vect{c}_0) + \mSighat_{l,l\vert l} \mSighat_{l,l\vert l-1}^{-1} \vmuhat_{l\vert l-1}.
\end{equation}
By the matrix inversion lemma, we obtain \cref{eqn:forward_pass_covariance_update}, that is,
\begin{equation}
    \mSighat_{l,l\vert l} 
    = \mSighat_{l,l\vert l-1} - \mSighat_{l,l\vert l-1} \mat{\tilde{C}}^T \big( \mSig_{\vv} + \mat{\tilde{C}} \mSighat_{l,l\vert l-1} \mat{\tilde{C}}^T \big)^{-1}  \mat{\tilde{C}} \mSighat_{l,l\vert l-1}.
\end{equation}
To compute the require quantities $\vmuhat_{l\vert l-1}$ and $\mSighat_{l,l\vert l-1}$, we use \cref{eqn:E_step_dynamical_model_app} and the fact that $\vw_{l-1}$ and $\vect{y}_{0:l-1}$ are independent, yielding
\begin{equation}
\begin{aligned}
    \vmuhat_{l\vert l-1} &= \mathbb{E}\big[ \vect{z}_l \vert \vect{y}_{0:l-1} \big] \\
    &= \mathbb{E}\big[ \mat{A}_{l-1} \vect{z}_{l-1} + \vect{b}_{l-1} + \vw_{l-1} \vert \vect{y}_{0:l-1} \big] \\
    &= \mat{A}_{l-1} \vmuhat_{l-1\vert l-1} + \vect{b}_{l-1}
\end{aligned}
\end{equation}
and
\begin{equation}
\begin{aligned}
    \mSighat_{l,l\vert l-1} &= \mathbb{E}\big[ (\vect{z}_l - \mat{A}_{l-1} \vmuhat_{l-1\vert l-1} - \vect{b}_{l-1}) (\vect{z}_l - \mat{A}_{l-1} \vmuhat_{l-1\vert l-1} - \vect{b}_{l-1})^T \vert \vect{y}_{0:l-1} \big] \\
    &= \mathbb{E}\big[ ( \mat{A}_{l-1} (\vect{z}_{l-1} - \vmuhat_{l-1\vert l-1}) + \vw_{l-1}) ( \mat{A}_{l-1} (\vect{z}_{l-1} - \vmuhat_{l-1\vert l-1}) + \vw_{l-1})^T \vert \vect{y}_{0:l-1} \big] \\
    &= \mat{A}_{l-1} \mSighat_{l-1,l-1\vert l-1} \mat{A}_{l-1}^T + \mSig_{\vw},
\end{aligned}
\end{equation}
which is \cref{eqn:forward_pass_covariance_increment}.
Combining the last four expressions, using the Kalman gain defined in \cref{eqn:Kalman_gain}, and collecting terms yields
\begin{multline}
    \vmuhat_{l\vert l} 
    = \mSighat_{l,l\vert l} \mat{\tilde{C}}^T \mSig_{\vv}^{-1} (\vect{y}_l - \vect{c}_0) 
    + \big( \mSighat_{l,l\vert l-1} - \mat{K}_l \mat{\tilde{C}} \mSighat_{l,l\vert l-1} \big) \mSighat_{l,l\vert l-1}^{-1} \vmuhat_{l\vert l-1} \\
    = \mSighat_{l,l\vert l-1} \mat{\tilde{C}}^T \left[ \mat{I} - \big( \mSig_{\vv} + \mat{\tilde{C}} \mSighat_{l,l\vert l-1} \mat{\tilde{C}}^T \big)^{-1}  \mat{\tilde{C}} \mSighat_{l,l\vert l-1} \mat{\tilde{C}}^T \right] \mSig_{\vv}^{-1} (\vect{y}_l - \vect{c}_0) 
    + \vmuhat_{l\vert l-1} 
    - \mat{K}_l \mat{\tilde{C}} \vmuhat_{l\vert l-1}.
\end{multline}
Using the identity $\big[ \mat{I} - (\mat{A} + \mat{B})^{-1} \mat{B} \big] \mat{A}^{-1} = (\mat{A} + \mat{B})^{-1}$ (see \cite{Yu2004derivation}) gives
\begin{equation}
    \vmuhat_{l\vert l}
    = \underbrace{\mSighat_{l,l\vert l-1} \mat{\tilde{C}}^T \big( \mSig_{\vv} + \mat{\tilde{C}} \mSighat_{l,l\vert l-1} \mat{\tilde{C}}^T \big)^{-1}}_{\mat{K}_l} (\vect{y}_l - \vect{c}_0) 
    + \vmuhat_{l\vert l-1} 
    - \mat{K}_l \mat{\tilde{C}} \vmuhat_{l\vert l-1},
\end{equation}
which, after collecting terms, is \cref{eqn:forward_pass_mean_update}.

To complete the derivation of the forward recursion, it remains to find $\vmuhat_{0,0}$ and $\mSighat_{0,0\vert 0}$.
Bayes' theorem gives
\begin{equation}
    \log P(\vect{z}_0 \vert \vect{y}_0) = \log P(\vect{y}_0 \vert \vect{z}_0) + \log P(\vect{z}_0) - \log P(\vect{y}_0).
\end{equation}
Proceeding as we did before and equating the linear and quadratic terms in $\vect{z}_0$ yields
\begin{equation}
    \mSighat_{0,0\vert 0}^{-1} = \mat{\tilde{C}}^T \mSig_{\vv}^{-1} \mat{\tilde{C}} + \mSig_0^{-1}
    \quad \mbox{and} \quad
    \vmuhat_{0\vert 0} = \mSighat_{0,0\vert 0} \mat{\tilde{C}}^T \mSig_{\vv}^{-1} (\vect{y}_0 - \vect{c}_0) + \mSighat_{0,0\vert 0} \mSig_0^{-1} \vmu_0.
\end{equation}
Applying the matrix inversion lemma gives
\begin{equation}
    \mSighat_{0,0\vert 0} 
    = \mSig_0 - \mSig_0 \mat{\tilde{C}}^T \big( \mSig_{\vv} + \mat{\tilde{C}} \mSig_0 \mat{\tilde{C}}^T \big)^{-1}  \mat{\tilde{C}} \mSig_0,
\end{equation}
which is equivalent to setting $\mSighat_{0,0\vert -1} = \mSig_0$ in \cref{eqn:forward_pass_covariance_update}.
Using this in the definition of the Kalman gain $\mat{K}_0$ and proceeding as above, we obtain
\begin{equation}
\begin{aligned}
    \vmuhat_{0\vert 0} 
    &= \mSighat_{0,0\vert 0} \mat{\tilde{C}}^T \mSig_{\vv}^{-1} (\vect{y}_0 - \vect{c}_0) 
    + \vmu_0 - \mat{K}_0  \mat{\tilde{C}} \vmu_0 \\
    &= \mSig_0 \mat{\tilde{C}}^T \left[ \mat{I} - \big( \mSig_{\vv} + \mat{\tilde{C}} \mSig_0 \mat{\tilde{C}}^T \big)^{-1}  \mat{\tilde{C}} \mSig_0 \mat{\tilde{C}}^T \right] \mSig_{\vv}^{-1} (\vect{y}_0 - \vect{c}_0) 
    + \vmu_0 - \mat{K}_0  \mat{\tilde{C}} \vmu_0 \\
    &= \mSig_0 \mat{\tilde{C}}^T \big( \mSig_{\vv} + \mat{\tilde{C}} \mSig_0 \mat{\tilde{C}}^T \big)^{-1} (\vect{y}_0 - \vect{c}_0) 
    + \vmu_0 - \mat{K}_0  \mat{\tilde{C}} \vmu_0 \\
    &= \vmu_0 + \mat{K}_0 (\vect{y}_0 - \vect{c}_0 - \mat{\tilde{C}} \vmu_0),
\end{aligned}
\end{equation}
which is equivalent to setting $\vmuhat_{0\vert -1} = \vmu_0$ in \cref{eqn:forward_pass_mean_update}.
This completes the derivation of the forward recursion.

\subsection{Backward recursion (smoothing)}
The desired quantities $\vmuhat_l$, $\mSighat_{l,l}$, and $\mSighat_{l, l+1}$ are parameters of the Gaussian probability distribution $P(\vect{z}_{l+1}, \vect{z}_1 \vert \vect{y}_{0:L})$.
For $l=L$, the quantities $\vmuhat_L = \vmuhat_{L\vert L}$ and $\mSighat_{L,L} = \mSighat_{L,L\vert L}$ were obtained during the forward recursion.
Beginning with $l=L-1$ we use the known $\vmuhat_{l+1}$ and $\mSighat_{l+1,l+1}$ to determine $\vmuhat_l$, $\mSighat_{l,l}$, and $\mSighat_{l, l+1}$ using a recursive process.
Following \cite{Yu2004derivation}, Bayes' theorem and the Markov property of \cref{eqn:E_step_dynamical_model_app} gives
\begin{equation}
\begin{aligned}
    \log P(\vect{z}_{l+1}, \vect{z}_l \vert \vect{y}_{0:L})
    &= \log P( \vect{z}_l \vert \vect{z}_{l+1}, \vect{y}_{0:L}) + \log P(\vect{z}_{l+1} \vert \vect{y}_{0:L}) \\
    &= \log P( \vect{z}_l \vert \vect{z}_{l+1}, \vect{y}_{0:l}) + \log P(\vect{z}_{l+1} \vert \vect{y}_{0:L}) \\
    &= \log P(\vect{z}_{l+1} \vert \vect{z}_l) + \log P( \vect{z}_l \vert \vect{y}_{0:l}) - \log P( \vect{z}_{l+1} \vert \vect{y}_{0:l}) + \log P(\vect{z}_{l+1} \vert \vect{y}_{0:L}).
\end{aligned}
\end{equation}
The parameters $\vmuhat_l$ and
\begin{equation}
    \begin{bmatrix}
        \mat{S}_{1,1} & \mat{S}_{1,2} \\
        \mat{S}_{2,1} & \mat{S}_{2,2}
    \end{bmatrix} = 
    \begin{bmatrix}
        \mSighat_{l+1,l+1} & \mSighat_{l+1,l} \\
        \mSighat_{l,l+1} & \mSighat_{l,l}
    \end{bmatrix}^{-1},
\end{equation}
are determined by the following terms of
\begin{multline}
    \log P(\vect{z}_{l+1}, \vect{z}_l \vert \vect{y}_{0:L}) 
    = -\frac{1}{2} 
    \begin{bmatrix} (\vect{z}_{l+1} - \vmuhat_{l+1})^T & 
    (\vect{z}_l - \vmuhat_l)^T \end{bmatrix}
    \begin{bmatrix}
        \mat{S}_{1,1} & \mat{S}_{1,2} \\
        \mat{S}_{2,1} & \mat{S}_{2,2}
    \end{bmatrix}
    \begin{bmatrix} \vect{z}_{l+1} - \vmuhat_{l+1} \\ 
    \vect{z}_l - \vmuhat_l \end{bmatrix} + \cdots \\
    = -\frac{1}{2}\vect{z}_{l+1}^T \mat{S}_{1,1} \vect{z}_{l+1}
    - \vect{z}_{l+1}^T \mat{S}_{1,2} \vect{z}_{l}
    -\frac{1}{2}\vect{z}_{l}^T \mat{S}_{2,2} \vect{z}_{l}
    + \vect{z}_{l}^T\big( \mat{S}_{2,1} \vmuhat_{l+1} + \mat{S}_{2,2} \vmuhat_{l} \big)
    + \cdots.
\end{multline}
We compute the corresponding terms of
\begin{multline}
    \log P(\vect{z}_{l+1} \vert \vect{z}_l) 
    = -\frac{1}{2}\big( \vect{z}_{l+1} - \mat{A}_l\vect{z}_l - \vect{b}_l \big)^T \mSig_{\vw}^{-1} \big( \vect{z}_{l+1} - \mat{A}_l\vect{z}_l - \vect{b}_l \big) + \cdots \\
    = -\frac{1}{2}\vect{z}_{l+1}^T \mSig_{\vw}^{-1} \vect{z}_{l+1} + \vect{z}_{l+1}^T \mSig_{\vw}^{-1} \mat{A}_l \vect{z}_l -\frac{1}{2} \vect{z}_l^T \mat{A}_l^T \mSig_{\vw}^{-1} \mat{A}_l \vect{z}_l - \vect{z}_l^T\mat{A}_l^T \mSig_{\vw}^{-1} \vect{b}_l  + \cdots
\end{multline}
\begin{equation}
    \log P( \vect{z}_l \vert \vect{y}_{0:l})
    = -\frac{1}{2} \vect{z}_l^T \mSighat_{l,l\vert l}^{-1} \vect{z}_l + \vect{z}_l^T \mSighat_{l,l\vert l}^{-1} \vmuhat_{l\vert l} + \cdots
\end{equation}
\begin{equation}
    \log P( \vect{z}_{l+1} \vert \vect{y}_{0:l})
    = -\frac{1}{2} \vect{z}_{l+1}^T \mSighat_{l+1,l+1\vert l}^{-1} \vect{z}_{l+1} + \cdots
\end{equation}
\begin{equation}
    \log P(\vect{z}_{l+1} \vert \vect{y}_{0:L})
    = -\frac{1}{2} \vect{z}_{l+1}^T \mSighat_{l+1,l+1}^{-1} \vect{z}_{l+1} + \cdots.
\end{equation}
Matching the corresponding terms in the above equations we obtain
\begin{equation}
    \begin{bmatrix}
        \mSighat_{l+1,l+1} & \mSighat_{l+1,l} \\
        \mSighat_{l,l+1} & \mSighat_{l,l}
    \end{bmatrix}^{-1}
    = 
    \begin{bmatrix}
        \mSig_{\vw}^{-1} - \mSighat_{l+1,l+1\vert l}^{-1} + \mSighat_{l+1,l+1}^{-1} & -\mSig_{\vw}^{-1} \mat{A}_l \\
        -\mat{A}_l^T \mSig_{\vw}^{-1} & \mat{A}_l^T \mSig_{\vw}^{-1} \mat{A}_l + \mSighat_{l,l\vert l}^{-1}
    \end{bmatrix}
    \label{eqn:backward_sigma_matching}
\end{equation}
and
\begin{equation}
    \mat{S}_{2,1} \vmuhat_{l+1} + \mat{S}_{2,2} \vmuhat_{l}
    = -\mat{A}_l^T \mSig_{\vw}^{-1} \vmuhat_{l+1} + \big( \mat{A}_l^T \mSig_{\vw}^{-1} \mat{A}_l + \mSighat_{l,l\vert l}^{-1} \big) \vmuhat_{l}
    = \mSighat_{l,l\vert l}^{-1} \vmuhat_{l\vert l} - \mat{A}_l^T \mSig_{\vw}^{-1} \vect{b}_l.
    \label{eqn:backward_mean_matching}
\end{equation}
Inversion of \cref{eqn:backward_sigma_matching} is carried out in \cite{Yu2004derivation}, yielding the recursion formulas \cref{eqn:backward_pass_covariance_update} and the expressions
\begin{equation}
    \mat{S}_{2,2}^{-1} = \mSighat_{l,l\vert l} - \mat{J}_l \mSighat_{l+1,l+1\vert l} \mat{J}_l^T
    \quad \mbox{and} \quad
    \mat{S}_{2,2}^{-1}\mat{S}_{2,1} = - \mat{J}_l,
\end{equation}
where $\mat{J}_l = \mSighat_{l,l\vert l} \mat{A}_l^T \mSighat_{l+1,l+1\vert l}^{-1}$.
Solving \cref{eqn:backward_mean_matching} for $\vmuhat_l$ yields
\begin{equation}
\begin{aligned}
    \vmuhat_l &= \mat{S}_{2,2}^{-1}\mSighat_{l,l\vert l}^{-1} \vmuhat_{l\vert l} - \mat{S}_{2,2}^{-1}\mat{S}_{2,1} (\vmuhat_{l+1}-\vect{b}_l) \\
    &= \vmuhat_{l\vert l} - \mat{J}_l \underbrace{\mSighat_{l+1,l+1\vert l} \mat{J}_l^T \mSighat_{l,l\vert l}^{-1}}_{\mat{A}_l} \vmuhat_{l\vert l} + \mat{J}_l (\vmuhat_{l+1}-\vect{b}_l),
\end{aligned}
\end{equation}
which is the recursion formula \cref{eqn:backward_pass_mean_update}.
This completes the derivation of the backward recursion.

\end{document}